\documentclass[11pt,a4paper]{article}
\usepackage[]{times}
\usepackage{fullpage}
\usepackage{graphicx}
\usepackage{float}
\usepackage{subfigure}
\usepackage{amsmath}
\usepackage{fancyhdr}
\usepackage[colorlinks]{hyperref}
\date{}

\newcommand{\prov}{{\sc Proof}.\hspace*{3mm} }

\newcommand{\QED}{$\rule{2mm}{2mm}$}

\newcommand{\natu}{{\sf I \! N}}
\newtheorem{theorem}{Theorem}[section]
\newtheorem{lemma}[theorem]{Lemma}
\newtheorem{e-proposition}[theorem]{Proposition}

\newtheorem{e-definition}[theorem]{Definition\rm}
\newtheorem{remark}{\it Remark\/}

\title{Optimal-rate Lagrange and Hermite finite elements for Dirichlet \\ problems in curved domains with straight-edged triangles}

\author{
    Vitoriano Ruas$^{1}$\thanks{Sorbonne Universit\'e, Campus de Jussieu, case 162, 4 place jussieu, Couloir 55-65, 4\`eme \'etage, 75005 Paris, France.}
		\\[1mm]
  {\small $^{1}$ Institut Jean Le Rond d'Alembert, CNRS UMR 7190, 
  Sorbonne Universit\'e, Paris, France.}\\[1mm]
  {\small e-mail: {\it vitoriano.ruas@upmc.fr}}}

\begin{document}
\maketitle

\begin{abstract}
One of the reasons for the success of the finite element method is its versatility to deal with different types of geometries. This is particularly true of problems posed in curved domains of arbitrary shape. In the case of second order boundary-value problems with Dirichlet conditions prescribed on curvilinear boundaries, method's isoparametric version for meshes consisting of curved triangles or tetrahedra has been mostly employed to recover the optimal approximation properties known to hold for methods of order greater than one based on standard straight-edged elements, in the case of polygonal or polyhedral domains. However, besides algebraic and geometric inconveniences, the isoparametric technique is limited in scope, since its extension to degrees of freedom other than function values is not straightforward. 
The purpose of this paper is to study a simple alternative that bypasses the above drawbacks, without eroding qualitative approximation properties. Besides other advantages, this technique can do without curved elements and is based only on polynomial algebra. It is first illustrated in the case of the convection-diffusion equation solved with standard Lagrange elements. Then it is applied to the solution with Hermite elements of the biharmonic equation with Dirichlet boundary conditions.\\

\noindent \textbf{Keywords:} Biharmonic; Curvilinear boundary; Dirichlet; Finite elements; Hermite; Lagrange; Straight-edged; Triangles.\\

\noindent \textbf{AMS Subject Classification: 65N30, 74S05, 76M10, 78M10, 80M10.}
\end{abstract}

\section{Study framework}

This work deals with a new method for solving boundary-value problem posed in a smooth curved two-dimensional domain of arbitrary shape. For the three-dimensional version of his method the author refers to 
\cite{arXiv3D} and \cite{IMAJNA}. \\
It is well known that in the finite-element solution of elliptic equations with Dirichlet conditions on a curvilinear boundary, a considerable order lowering may occur if prescribed boundary values are shifted to nodes that are not mesh vertexes of an approximating polygon or polyhedron 
formed by the union of straight-edged $N$-simplexes of a fitted mesh. Over four decades ago some techniques were designed in order to remedy such a loss of accuracy, and possibly attain the same theoretical optimal orders as in the case of a polytopic domain, assuming that the solution is sufficiently smooth. Two examples of such attempts are the \textit{interpolated boundary condition method} by Nitsche and Scott (cf. \cite{Nitsche} and \cite{Scott}), and the method introduced by Zl\'amal in \cite{Zlamal} and extended by \v{Z}\'eni\v{s}ek in \cite{Zenisek}.\\
The principle our method is based upon is close to the interpolated boundary conditions 
studied in \cite{BrennerScott} for two-dimensional problems. Although the latter technique is very intuitive and is known since the seventies (cf. 
\cite{Scott}), it has been of limited use so far. Among the reasons for this we could quote its difficult implementation, the lack of an extension to three-dimensional problems, and most of all, restrictions on the choice of boundary nodal points to reach optimal convergence rates. In contrast our method is simple to implement in both two- and three-dimensional geometries. Moreover optimality is attained very naturally in both cases for various choices of boundary nodal points, as seen hereafter. \\ 
Since long the isoparametric version of the finite element method for meshes consisting of curved triangles or tetrahedra (cf. \cite{Zienkiewicz}), has been considered as the ideal way to handle Dirichlet conditions prescribed on curved boundaries. It turns out that, besides a more elaborated description of the mesh, the isoparametric technique inevitably leads to the integration of rational functions to compute the system 
matrix. This raises the delicate question on how to choose the right numerical quadrature formula in the master element in the case of complex non linear problems. In contrast, in the technique studied in this paper exact numerical integration can be used for this purpose, at least in the most frequent situations, since we only have to deal with polynomial integrands. Moreover the element geometry remains the same as in the case of polygonal or polyhedral domains. It is noteworthy that both advantages are conjugated with the fact that no erosion of qualitative approximation properties results from the application of our technique, as compared to the equivalent isoparametric one. 
We should also emphasize that this approach is particularly handy, whenever the finite element method under consideration has normal components or normal derivatives as degrees of freedom. Indeed in this case isoparametric analogs are either not so easy to define (see. e.g. \cite{RT}) or are simply unknown. \\

\indent An outline of the paper is as follows. Section 2 is first devoted to the model problem in a smooth two-dimensional domain 
selected for the presentation of our method, namely, the convection-diffusion equation. Some pertaining notations are also given therein, followed by a preliminary material concerning the boundary of this domain, as related to the family of meshes considered in the sequel.
In Section 3 we present our method as applied to solve the model problem  
with Dirichlet boundary conditions; higher order conforming Lagrange finite elements based on meshes with straight-edged triangles are studied and corresponding well-posedness results are demonstrated. In Section 4 we prove error estimates in the $H^1$-norm for the method introduced in the previous section. Moreover $L^2$-error estimates are demonstrated in relevant cases, which to the best of author's knowledge are unprecedented for the class of problems addressed in this work. 
In Section 5 we show that the principles presented in Section 3 extend very naturally to fourth-order problems solved by Hermite finite element methods with normal-derivative degrees of freedom.
We conclude in Section 6 with some comments on the whole work. \\
Numerics is not the focus of this work. For numerical experimentation the author refers to several other papers of his quoted in the bibliography. Nevertheless, with the aim of dissipating any skepticism 
about the performance of our method vis-\`a-vis classical techniques, some comparative results 
are supplied in the Appendix. 

\section{Preliminaries}

The methodology to enforce Dirichlet boundary conditions on curvilinear boundaries studied in this work applies to many types of equations. However, in order to avoid non essential difficulties, we consider as a model the following convection-diffusion equation in a two-dimensional smooth domain $\Omega$ with boundary $\Gamma$, namely:
\begin{equation}
\label{Poisson}
\left\{
\begin{array}{l}
 -\nu \Delta u + {\bf b} \cdot {\bf grad} \; u = f \mbox{ in } \Omega \\
 u = g \mbox{ on } \Gamma,
\end{array}
\right.
\end{equation} 
\noindent where $\nu$ is the diffusion coefficient and ${\bf b}$ is a given continuous and solenoidal velocity field defined in $\Re^2$. $f$ and $g$ in turn are given functions defined in $\Omega$ and on $\Gamma$, having suitable regularity properties. We shall be dealing with approximation 
methods of order $k$ for $k > 1$ in the standard (semi-)norm $\parallel {\bf grad} (\cdot) \parallel_{0}$ of $H^1$, as long as $u \in H^{k+1}(\Omega)$, where 
$\parallel \cdot \parallel_{0}$ denotes the standard norm of $L^2(\Omega)$.
Accordingly, we shall assume that $f \in H^{k-1}(\Omega)$ and $g \in H^{k+1/2}(\Gamma)$ (cf. \cite{Adams}). Although the method to be described below applies to any 
$g$, for the sake of simplicity henceforth we shall take $g \equiv 0$. In this case, for the assumed regularity of $u$ to hold, we require that both ${\bf b}$ and $\Gamma$ be sufficiently smooth and at least of the $C^{k-1}$-class. 
Actually, more than this, if $k=2$ we make the assumption that the curvature of $\Gamma$ (cf. \cite{Cartan}) is uniquely defined almost everywhere. Eventually, for $k>2$ too we will require more from $\Gamma$ than being of the $C^{k-1}$-class.

\subsection{Meshes and related sets} 

Let us be given a mesh ${\mathcal T}_h$ consisting of straight-edged triangles satisfying the usual compatibility conditions (see e.g. \cite{Ciarlet}). Every element of ${\mathcal T}_h$ is to be viewed as a closed set. Moreover this mesh is assumed to fit $\Omega$ in such a way that all the vertexes of the polygon $\cup_{T \in {\mathcal T}_h} T$ lie on $\Gamma$. We denote the interior of this union set by $\Omega_h$ and define $\tilde{\Omega}_h := \Omega \cup \Omega_h$ and $\Omega^{'}_h := \Omega \cap \Omega_h$. The boundaries of $\Omega_h$ and $\tilde{\Omega}_h$ are respectively denoted by $\Gamma_h$ and $\tilde{\Gamma}_h$ and moreover we set   
$\Gamma^{'}_h:= \Omega_h \cap \Gamma$. ${\mathcal T}_h$ is assumed to belong to a regular family of partitions 
in the sense of \cite{Ciarlet} (cf. Section 3.1), though not necessarily quasi-uniform. \\
The diameter of every $\forall T \in {\mathcal T}_h$ is represented by $h_T$,  while $h :=  \max_{T \in {\mathcal T}_h} h_T$. We make the non essential and yet reasonable assumption that any element in ${\mathcal T}_h$ have at most one edge  contained in $\Gamma_h$. Actually such a condition is commonly fulfilled in practice, so that excessively flat triangles are avoided. \\
Let ${\mathcal S}_h$ be the subset of ${\mathcal T}_h$ consisting of triangles $T$ having 
one edge on $\Gamma_h$, say $e_T$. For every $T \in {\mathcal S}_h$ we denote by $O_T$ the vertex of $T$ not belonging to $\Gamma$; moreover we define $T^{\Gamma}$ to be the curved triangle delimited by $\Gamma$ and the two edges of $T$ intersecting at $O_T$. 
Notice that, owing to our initial assumption, no triangle in ${\mathcal T}_h \setminus {\mathcal S}_h$ has a nonempty intersection with $\Gamma_h$. We denote by ${\mathcal Q}_h$ the subset of ${\mathcal S}_h$ consisting of elements $T$ 
such that $T \setminus \Omega$ is not restricted to a pair of vertexes of $\Gamma_h$.

\subsection{Notations} 

Hereafter $\parallel \cdot \parallel_{r,D}$ and $| \cdot |_{r,D}$ represent, respectively, the standard norm and semi-norm of Sobolev space $H^{r}(D)$ (cf. \cite{Adams}), for  
$r \in \Re^{+}$ with $H^0(D)=L^2(D)$, $D$ being a subset of the closure of $\tilde{\Omega}_h$. We also denote by $\parallel \cdot \parallel_{m,p,D}$ and $| \cdot |_{m,p,D}$ the usual norm and semi-norm of $W^{m,p}(D)$ for $m \in \natu^{*}$ and 
$p \in [1,\infty] \setminus \{2\}$ with $W^{0,p}(D)=L^p(D)$, and also for $W^{m,2}(D)=H^m(D)$, whenever convenient. In case $D=\Omega$ the subscript $D$ is dropped. Throughout this article ${\mathcal P}_k(D)$ represents the space of polynomials of degree less than or equal to $k$ defined in $D$.\\
Henceforth we denote by $D^j w$ the $j$-th order tensor whose components are the $j$-th order partial derivatives with respect to the space variables of a function $w$ in the strong or the weak sense. Alternatively we may also write $H(w)$ instead of $D^2 w$ and ${\bf grad}\;w$ instead of $D^1 w$. \\
Finally we introduce the notations $\parallel \cdot \parallel_{0,h}$ and $\parallel \cdot \parallel^{'}_{0,h}$ 
for the standard norms of $L^2(\Omega_h)$ and $L^2(\Omega^{'}_h)$, respectively, 
 which will play a key role in the reliability analysis of our method. This is because all our error estimates will be given in the former norm if $\Omega$ is convex and in the latter otherwise. \\
In this respect it is noticeable that for a given mesh and a function $v \in L^2(\Omega)$, $\parallel v \parallel_{0,h}$ (resp. $\parallel \cdot \parallel^{'}_{0,h}$) may equal zero, even if $v$ does not vanish in 
$\Omega \setminus \Omega_h$. However in asymptotic terms this situation is ruled out as far as $u$ is concerned. Indeed the estimates are supposed to hold as $h$ goes to zero, since the family of meshes under consideration is regular (cf. \cite{Ciarlet}, Sect. 3.1). Thus the meshes asymptotically cover the whole $\Omega$. Incidentally this apparently 
indefinite error measure in the case of curved domains is the one used in classical textbooks on the mathematical analysis of the finite element method, such as \cite{Ciarlet} (cf. Section 4.4. p.266 and on) and \cite{StrangFix} (cf. Section 4.4, p.192 and on).

\subsection{Basic assumptions for the formal analysis}
 
Although this is by no means necessary, neither to define our method, nor to implement it, 
henceforth we assume that the meshes in use are fine enough to satisfy some geometric criteria. This assumption is a key sufficient condition for the subsequent reliability results to hold. It also enables the capture of all the nuances of $\Gamma$ by its discrete counterpart $\Gamma_h$, taking advantage of the great flexibility of triangular meshes to fit curvilinear boundaries, even those with sharp variations of shape. \\

We first require the following condition:\\

\noindent \underline{\textit{Assumption}$^{+}$ :} Let $T_{\rho}$ be the homothetic transformation of $T$ with center $O_T$ and ratio $\rho <1$, possibly small. $h$ is small enough for the intersection $P$ with $\Gamma$ belonging to $T^{\Gamma}$ of a straight line joining any point of $T_{\rho}$ to a point $M \in e_T$ to be uniquely defined $\forall T \in {\mathcal S}_h$.  \rule{2mm}{2mm} \\

In addition to \textit{Assumption}$^{+}$ the following condition is also supposed to be satisfied by the meshes: \\
Let $Q_T$  be the closest intersection with $\Gamma$ of the perpendicular to $e_T$ passing through its mid-point $M_T$. We know that there exists a ball $B(Q_T,r_T)$ and a straight line $\Pi_T$ swept by the coordinate $x_T$ of an orthogonal coordinate system $(O,x_T,y_T)$ with a suitably chosen origin $O$, such that a function $f_T(x_T)$ of the piecewise $C^2$-class uniquely expresses the coordinate $y_T$ of points located on $\Gamma$, as long as they lie in $B(Q_T,r_T)$ (cf. \cite{Evans}).\\ 

\noindent \underline{\textit{Assumption}$^{*}$ :} $h$ is small enough for $\Pi_T$ to be aligned with $e_T$ and the ball $B(Q_T,r_T)$ to contain $e_T$ $\forall T \in {\mathcal S}_h$  \rule{2mm}{2mm} \\

Some important consequences of both assumptions above are as follows:

\begin{e-proposition}
\label{prop01}
If \textit{Assumption}$^{+}$ and \textit{Assumption}$^{*}$ hold there exists a constant $C_{\Gamma}$ depending only on $\Gamma$ such that 
$\forall M \in e_T$ the length of the segment joining $M$ and $P \in T^{\Gamma} \cap \Gamma$ aligned with $O_T$ and $M$ is bounded above by $C_{\Gamma} h_T^2$.
\end{e-proposition}

\prov Denoting by $l_T$ the length of $e_T$,  
let $x$ be the abscissa along $e_T$ in the interval $[0,l_T]$, whose orientation plays no role. From 
\textit{Assumption}$^{*}$ the portion of $\Gamma$ comprised between the two vertexes of $T$ lying on $\Gamma$ can be uniquely 
represented by a function $f_T$ of $x$, such that any point $P$ of such a portion of $\Gamma$ has coordinates $(x,f_T(x))$ in the cartesian coordinate system $(O,x,y)$, whose origin $O \in \Gamma$ is one of the end-points of $e_T$.\\
Next we prove that $|f_T^{''}|$ is uniformly bounded in $[0,l_T]$ independently of both $T \in {\mathcal T}_h$ and the mesh size $h$. With this aim we first recall that the  curvature $\kappa$ of $\Gamma$ at a point $P \in \Gamma$ with coordinates $(x,f_T(x))$ can be locally expressed in terms of $f_T$, in such a way that (see e.g. \cite{Goldman}):
\begin{equation}
\label{Label} 
|\kappa(P)| = \displaystyle \frac{|f_T^{''}(x)|}{[1+|f_T^{'}(x)|^2]^{3/2}} \; \forall x \in [0,l_T]. 
\end{equation}
The proof is based on the fact that, provided $h$ is sufficiently small, the $L^{\infty}(0,l_T)$-norm of the function $f_T^{''}$ is bounded above by an expression depending only on the maximum absolute value of the curvature of $\Gamma$, multiplied by a constant independent of $T$. Next we give a rigorous justification of this assertion. \\         
Let ${\mathcal C}_{max} := \max_{P \in \Gamma} |\kappa(P)|$ and ${\mathcal H}_{max}:=\max_{x \in [0,l_T]} |f_T^{''}(x)|$. Since $f_T(0)=f_T(l_T)=0$, there is necessarily an abscissa $x_0 \in [0,l_T]$ at which $f_T^{'}$ vanishes, and hence we can write $|f^{'}_T(x)| = |\int_{x_0}^x f^{''}_T(s) ds|$ for $x \in [0,l_T]$. Then from \eqref{Label} straightforward calculations yield,

\[ |f_T^{''}(x)|^2 \leq {\mathcal C}_{max}^2 (1+l_T^2 {\mathcal H}_{max}^2)^{3} \; \forall x \in [0,l_T].\]

Now we assume that $l_T \leq \sqrt{\beta} / {\mathcal C}_{max}$, where $\beta$ is less than or equal to $4/27$.  
This means that the upper bound for ${\mathcal H}_{max}$ we are searching for satisfies,
\[ \displaystyle \frac{{\mathcal H}_{max}^2}{{\mathcal C}_{max}^2} \leq 1 + 3 \beta \displaystyle \frac{{\mathcal H}_{max}^2}{{\mathcal C}_{max}^2}+3 \beta^2 \displaystyle \frac{{\mathcal H}_{max}^4}{{\mathcal C}_{max}^4}+\beta^3 \displaystyle \frac{{\mathcal H}_{max}^6}{{\mathcal C}_{max}^6}.\]
For convenience we set $t:={\mathcal H}_{max}^2/{\mathcal C}_{max}^2$. Next we check whether there exists $t_1 > 0$ such that $0 \leq \varphi(t):= t(1-3\beta-3 \beta^2 t-\beta^3 t^2) \leq 1$ for every $t$ in $[0,t_1]$. Straightforward calculations show that if $\beta \leq 4/27$, the function $\varphi(t)$ is non negative for $0 \leq t \leq t_{max} := (-3 + \sqrt{4 \beta^{-1}-3})/(2 \beta)$ with $t_{max} > 0$.  
Moreover in the interval $[0,t_{max}]$ $\varphi$ attains a minimum at both $t=0$ and $t=t_{max}$ and only a local maximum greater than one at the point $t_0=(-3+\sqrt{\beta^{-1}})/(3 \beta)$. It follows that there exists a point $t_1 \in (0, t_{0})$ depending only on $\beta$ such that $\varphi(t_1) = 1$ and hence ${\mathcal H}_{max} \leq \sqrt{t_1} {\mathcal C}_{max}$. Notice that this upper bound is uniform and holds for all $T$ having an edge on $\Gamma_h$. \\
Let $N_0 \in \Gamma$ be the point of $e_T$ with abscissa $x_0$. Then at any point $N \in e_T$ we have $f_T^{'}(N) \leq {\mathcal H}_{max} length(\overline{N_0N}) \leq {\mathcal H}_{max} l_T$, which implies that 
${\mathcal G}_{max} := \parallel f_T^{'} \parallel_{0,\infty,e_T} \leq {\mathcal H}_{max} h_T$.  
Then again, the fact that $f_T(O)=0$ implies that $\forall N \in e_T$, $|f_T(N)| \leq length(\overline{ON}) {\mathcal G}_{max} $. Therefore $\forall N \in e_T$, $|f_T(N)| \leq {\mathcal G}_{max} h_T$ and thus $|f_T(N)| \leq {\mathcal H}_{max} h_T^2, \; \forall N \in e_T$.\\
At this point we define $\theta_{min}$ to be the minimum of the smallest angle of $T \in {\mathcal T}_h$ over the meshes in the regular family in use. A simple geometric argument allows us to conclude that the length of $\overline{MP}$ is bounded above by $f_T(M)/sin \theta_{min}$. Thus the result holds with $C_{\Gamma} = {\mathcal H}_{max} [sin \theta_{min}]^{-1}$.  \QED 

\begin{e-proposition}
\label{prop02}
Assume that $\Gamma$ is of the piecewise $C^{k+1}-class$ for $k>1$. Let $v^{(j)}$ denote the derivative of order $j$ with respect to $x$ of a sufficiently differentiable function $v(x)$, $k+1 \geq j \geq 0$ with $ v^{(0)} = v$, $ v^{(1)} = v^{'}$ and $v^{(2)} = v^{''}$.
If \textit{Assumption}$^{*}$ holds, there exist constants $C_{\Gamma}^{j}$ depending only of $\Gamma$ such that $|[f_T^{(j)}](M)| \leq 
C^j_{\Gamma} h_T^{\max[2-j,0]}$ $\forall M \in e_T$ for $j=0,1,2\ldots,k+1$.
\end{e-proposition} 

\prov From the proof of Proposition \ref{prop01} we infer that the result holds true with $C^j_{\Gamma} = {\mathcal H}_{max}$, for $j=0,1,2$.  Finally for $2 < j \leq k+1$ the bound is a simple consequence of the regularity assumptions on 
$\Gamma$. \rule{2mm}{2mm}

\section{Method description}

First of all we need some additional definitions regarding the skin $(\Omega \setminus \Omega_h) \cup (\Omega_h \setminus \Omega)$. 
$\forall T \in {\mathcal S}_h$ we denote by $\Delta_T$ the closed set delimited by $\Gamma$ and the edge $e_T$ of $T$ whose end-points belong to $\Gamma$, as illustrated in Figure 1.  
In this manner we can assert that, if $\Omega$ is convex, $\Omega_h$ is a proper subset of $\Omega$ and $\bar{\Omega}$ is the union of the disjoint sets $\Omega_h$ and $\displaystyle \cup_{T \in {\mathcal S}_h} \Delta_T$. Otherwise $\Omega_h \setminus \Omega$ is a nonempty set containing subsets of $T \in {\mathcal S}_h$ whose area is an $O(h_T^3)$ corresponding to non-convex portions of $\Gamma$. Whatever the case, the above configurations are of  merely academic interest and carry no practical meaning, as much as the sets 
$T^{\Delta}:=T \cup \Delta_T$ and $T^{'}:= T \cap \Omega$ $\forall T \in {\mathcal S}_h$.\\
Although this does not really play any role in practice, in order to avoid non essential technicalities, we will ideally consider that the mesh is constructed in such a way that convex and concave portions of $\Gamma$ correspond to convex and concave portions of $\Gamma_h$. This property is guaranteed if the mesh is sufficiently fine and moreover points separating such portions of $\Gamma$ are vertexes of polygon $\Omega_h$. In doing so, schematically if $e_T$ lies on a convex portion of $\Gamma_h$, $T$ is a proper subset of $T^{\Delta}$, while $T^{'}$ is a proper subset of $T$ if $e_T$ lies in a concave portion of $\Gamma_h$. \\ 

\begin{figure}[pb]
\label{fig1}
\centerline{\includegraphics[width=3.8in]{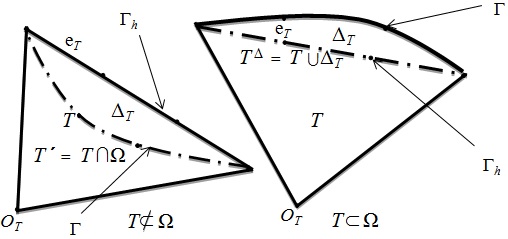}}
\vspace*{8pt}
\caption{Skin $\Delta_T$ related to a triangle $T$ next to a convex (right) or a concave (left) portion of $\Gamma$}
\end{figure} 

\indent Next we introduce two spaces $V_h$ and $W_h$ associated with ${\mathcal T}_h$. \\

$V_h$ is the standard Lagrange finite element space consisting of 
continuous functions $v$ defined in $\Omega_h$ that vanish on $\Gamma_h$, whose restriction to every $T \in {\mathcal T}_h$ is a polynomial of degree 
less than or equal to $k$ for $k \geq 2$. For convenience we extend every function $v \in V_h$ by zero in $\Omega \setminus \Omega_h$ . \\

$W_h$ is the space 
of functions defined in $\Omega_h$ having the properties listed below. \\
 
\begin{enumerate} 
\item The restriction of $w \in W_h$ to every $T \in {\mathcal T}_h$ is a polynomial of degree less than or equal to $k$;
\item Every $w \in W_h$ is continuous in $\Omega_h$ and vanishes at the vertexes of $\Gamma_h$; 
\item A function $w \in W_h$ is extended to $\Omega \setminus \Omega_h$ in such a way that its polynomial expression in $T \in {\mathcal S}_h$ also applies 
to points in $\Delta_T$;
\item $\forall T \in {\mathcal S}_h$, $w(P) = 0$ for every $P$ among the $k-1$ nearest intersections with $\Gamma$ 
of the line passing through the vertex $O_T$ of $T$ not belonging to $\Gamma_h$ and the points $M$ different from vertexes of $T$ subdividing the edge $e_T$ opposite to $O_T$ into $k$ segments of equal length (cf. Figure 2).
\end{enumerate}

\begin{figure}[pb]
\label{fig2}
\centerline{\includegraphics[width=3.8in]{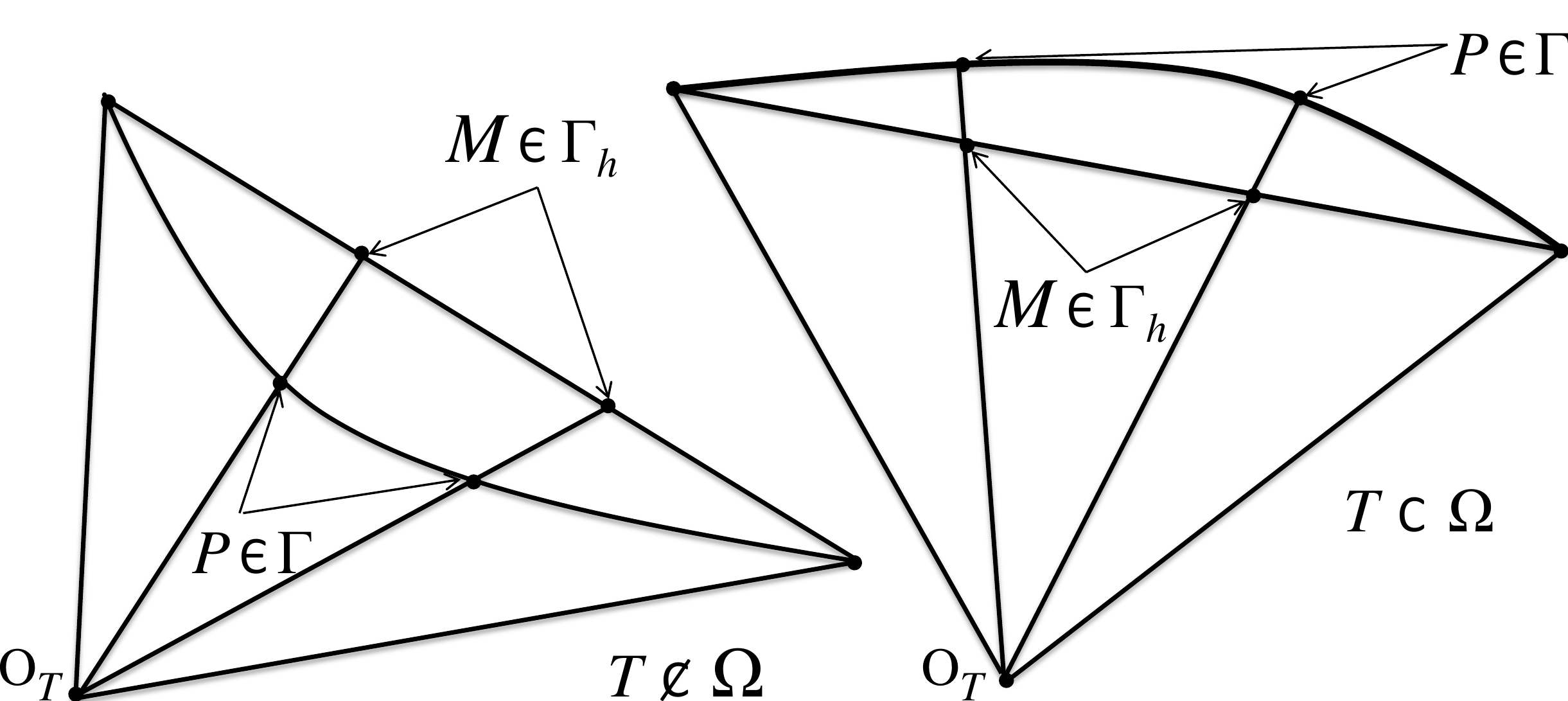}}
\vspace*{8pt}
\caption{Construction of nodes $P \in \Gamma$ for space $W_h$ related to Lagrangian nodes $M \in \Gamma_h$ for $k=3$}
\end{figure}
  
\begin{remark}
The construction of the nodes associated with $W_h$ located on $\Gamma$ advocated in item 4 is not mandatory. Notice that it differs 
from the intuitive construction of such nodes lying on normals to edges of $\Gamma_h$. The main advantage of this proposal is an easy determination of boundary-node coordinates by linearity, using a supposedly available analytical expression of $\Gamma$. 
Nonetheless the choice of boundary nodes ensuring our method's optimality is really wide.  \QED
\end{remark}

Henceforth the $n_k$ points of a triangle $T$ in ${\mathcal T}_h$ resulting from its subdivision into $k^2$ equal triangles will be called the lattice points of $T$ (cf. \cite{Zienkiewicz}). Notice that $n_k = (k+2)(k+1)/2$. \\

The fact that $W_h$ is a non empty finite-dimensional space is the consequence of the two following lemmata.

\begin{lemma}
\label{TDelta}
Provided $h$ satisfies \textit{Assumption}$^{+}$ and \textit{Assumption}$^{*}$ there exist two mesh-independent constants ${\mathcal C}_{\infty}$ and ${\mathcal C}_J$ depending only on $\Gamma$ and the shape regularity of ${\mathcal T}_h$ (cf. \cite{BrennerScott}, Ch.4, Sect. 4) such that  
$\forall w \in {\mathcal P}_k(T^{\Delta})$ and $\forall T \in {\mathcal S}_h$ it holds:  
\begin{equation}
\label{LinftyTDelta}
\parallel  w \parallel_{0,\infty,T^{\Delta}} \leq {\mathcal C}_{\infty} \parallel  w \parallel_{0,\infty,T^{'}} 
\end{equation}
and
\begin{equation}
\label{L2TDelta}
\parallel  w \parallel_{0,\infty,T^{\Delta}} \leq {\mathcal C}_J h_T^{-1} \parallel  w \parallel_{0,T^{'}}.
\end{equation}
\end{lemma}

\prov First we recall that the dimension of ${\mathcal P}_k(D)$ for any bounded open  
set $D$ of $\Re^2$ is $n_k$. \\
Let $0 < \lambda \leq 1$ be the largest possible value for the homothetic transformations $T^{'}_{\lambda}$ and $T^{\Delta}_{\lambda}$ of $T \in {\mathcal S}_h$ centered at $O_T$ and with ratios $\lambda$ and $\lambda^{-1}$, to be contained in $T^{'}$ and contain $T^{\Delta}$, respectively. 
Now set $\kappa^{'} :=1 - \sigma_{\mathcal T} C_{\Gamma} h_0$ and $\kappa^{\Delta}:= 1 + \sigma_{\mathcal T} C_{\Gamma} h_0$ as two numbers depending only on $\Gamma$, where $h_0$ is the largest value of $h$ such that \textit{Assumption}$^{+}$ and \textit{Assumption}$^{*}$ hold and $\kappa^{'}$ is not less than a certain number in the interval $(0,1]$, say $1/2$; $\sigma_{\mathcal T}$ in turn is a shape-regularity parameter of the family of meshes in use satisfying for every ${\mathcal T}_h$, $\sigma_{\mathcal T} \geq  
\max_{T \in {\mathcal T}_h} h_T/\eta_T$, where $\eta_T$ is the minimum height of $T$. From Proposition \ref{prop01}  
together with Thales' Proportionality Theorem, it is rather easy to infer that $\kappa^{'}$ and $\kappa^{\Delta}$ are such that the maximum diameters of triangles $T^{'}_{\lambda}$ and $T^{\Delta}_{\lambda}$ lie in the intervals $[\kappa^{'} h_T, h_T]$ and $[h_T,\kappa^{\Delta} h_T]$, respectively. 
Since both $T^{'}_{\lambda}$ and $T^{\Delta}_{\lambda}$ are similar 
to $T$, both triangles have the same shape regularity property as any other element in ${\mathcal T}_h$, provided the maximum diameter of each member of the family of triangulations in use is adjusted to take into account the thus modified maximum diameters.\\
Let $T \in {\mathcal S}_h$. Denoting by $\varphi_i$ the canonical basis function associated with the $i$-th lattice point $M_i$ of $T$ extended to $T^{\Delta}_{\lambda}$, for every $w \in {\mathcal P}_{k}(T^{\Delta})$ we can write,
\begin{equation}
\label{auxiliary6mono}
\parallel w \parallel_{0,\infty,T^{\Delta}} \leq \displaystyle \sum_{i=1}^{n_k} |w(M_i)| 
\max_{{\bf x} \in T_{\lambda}^{\Delta}} |\varphi_i({\bf x})|.
\end{equation}
Next we resort to the master triangle $\hat{T}$ with vertexes $(0,0),(1,0),(0,1)$ 
in a reference frame $(\hat{O},\hat{\bf x})$, where $\hat{O}$ corresponds to $O$. ${\mathcal F}_T$ being the affine mapping from $T$ onto $\hat{T}$, let $\hat{\varphi}_i$ and $\hat{w}$ be the transformations of $\varphi_i$ and $w$ under ${\mathcal F}_T$. Let also $\hat{T}^{'}_{\lambda}$ and $\hat{T}^{\Delta}_{\lambda}$ be the transformations of $T^{'}_{\lambda}$ and $T_{\lambda}^{\Delta}$ under ${\mathcal F}_T$. Then it holds: 
\begin{equation}
\label{auxiliary6bis}
\parallel w \parallel_{0,\infty,T^{\Delta}} \leq \hat{C}_1 \displaystyle \sum_{i=1}^{n_k} |w(M_i)| \; \forall w \in {\mathcal P}_{k}(T^{\Delta}),
\end{equation}
where $\hat{C}_1 =\displaystyle \max_{1 \leq i \leq n_k} \left[\max_{\hat{\bf x} \in \hat{T}_{\lambda}^{\Delta}} |\hat{\varphi }_i(\hat{\bf x})| \right].$\\
Owing to the equivalence of norms in the $n_k$-dimensional space ${\mathcal P}_k(\hat{T}^{'}_{\lambda})$, there exists a constant $\hat{C}_2$ depending only on $\hat{T}$, $\lambda$ and $k$ such that $\forall w \in {\mathcal P}_k(T^{\Delta})$, 
\begin{equation}
\label{auxiliary6ter}
\displaystyle \sum_{i=1}^{n_k} |w(M_i)| = \displaystyle \sum_{i=1}^{n_k} |\hat{w}({\mathcal F}_T(M_i))| \leq \hat{C}_2 \parallel \hat{w} \parallel_{0,\infty,\hat{T}^{'}_{\lambda}}.  
\end{equation} 
Combining \eqref{auxiliary6bis} and \eqref{auxiliary6ter} it easily follows that \eqref{LinftyTDelta} holds with 
${\mathcal C}_{\infty} = \hat{C}_1 \hat{C}_2$. \\
Finally we note that $area(\hat{T}^{'}_{\lambda}) \leq \hat{\mathcal C}_{J}^2 h_T^{-2} area(T^{'}_\lambda)$ with a constant $\hat{\mathcal C}_{J}$ independent of $T$. Then using again the equivalence of norms in ${\mathcal P}_k(\hat{T}^{'}_{\lambda})$ we infer the existence of another constant $\hat{C}_{\lambda}$ independent of $T$ for which it holds, 
\begin{equation}
\label{inverseTlambda}
\parallel \hat{w} \parallel_{0,\infty,\hat{T}^{'}_{\lambda}} \leq \hat{C}_{\lambda} 
\parallel \hat{w} \parallel_{0,\hat{T}^{'}_{\lambda}} \leq \hat{C}_{\lambda} \hat{\mathcal C}_{J} h_T^{-1} \parallel w \parallel_{0,T_{\lambda}^{'}} \; \forall w \in {\mathcal P}_k(T^{\Delta}).
\end{equation} 
Since $T^{'}_{\lambda} \subset T^{'}$, combining \eqref{auxiliary6bis}, \eqref{auxiliary6ter}, \eqref{inverseTlambda}, \eqref{L2TDelta} must hold with ${\mathcal C}_J= \hat{C}_{\lambda} \hat{\mathcal C}_{J} {\mathcal C}_{\infty}$ independently of $T^{'}$ and $T^{\Delta}$. \QED 

\begin{lemma}
\label{lemma1}
Provided $h$ is small enough, given a set of $m_k$ real values $b_{i}$, $i=1,\ldots,m_k$ with $m_k=(k+1)k/2$, $\forall T \in {\mathcal S}_h$ 
there exists a unique function $w_T \in {\mathcal P}_k(T)$ that vanishes at both 
vertexes of $T$ located on $\Gamma$ and at the $k-1$ points $P$ of $\Gamma$ defined in accordance with item 4. of the above definition of $W_h$, and takes the value $b_i$ respectively at the $m_k$ lattice points of $T$ not located on $\Gamma_h$. 
\end{lemma}   
 
\prov Let us first extend the vector $\vec{b}:=[b_1,b_2,\ldots,b_{m_k}]$ of $\Re^{m_k}$ into a vector of $\Re^{n_k}$ still denoted by $\vec{b}$, by adding $n_k-m_k=k+1$ zero components. If the boundary nodes $P$ were replaced by the corresponding $M \in \Gamma_h \cap T$, it is clear that the result would hold true, according to the well-known properties of Lagrange finite elements. The vector $\vec{a}$ of coefficients $a_i$ for $i=1,2,\ldots,n_k$ of the canonical basis functions 
$\varphi_i$ of ${\mathcal P}_k(T)$ for $1 \leq i \leq n_k$ would be given by $a_i=b_i$ for $1 \leq i \leq n_k$. Denoting by $M_i$ the lattice points of $T$, $i=1,2,\ldots,n_k$, this means that the matrix $K$ whose entries are $k_{ij} := \varphi_j(M_i)$ is the identity matrix. Let $\tilde{M}_i=M_i$ if $M_i \notin \Gamma \setminus \Gamma_h$ and $\tilde{M}_i$ be the node of the type $P$ associated with $M_i$ otherwise. The Lemma will be proved if the $n_k \times n_k$ linear system $\tilde{K} \vec{a} = \vec{b}$ is uniquely solvable, where $\tilde{K}$ is the matrix with entries $\tilde{k}_{ij}:=\varphi_j(\tilde{M}_{i})$. Clearly we have $\tilde{K} = K + 
E_K$, where the entries of $E_K$ are $e_{ij}:= \varphi_j(\tilde{M}_{i}) - \varphi_j(M_{i})$. \\
From Proposition \ref{prop01} the length of the segment $\overline{M_i\tilde{M}_i}$ is bounded above by $C_{\Gamma} h_T^2$. 
It follows from Rolle's Theorem that $\forall \;i,j$, 
$|e_{ij}| \leq C_{\Gamma} h_T ^2 \max_{{\bf x} \in T_{\Delta}} |{\bf grad}\;\varphi_j({\bf x})|$, or yet, recalling \eqref{LinftyTDelta},  
\begin{equation}
\label{boundeij}
|e_{ij}| \leq {\mathcal C}_{\infty} C_{\Gamma} h_T ^2 \max_{{\bf x} \in T} |{\bf grad}\;\varphi_j({\bf x})|.
\end{equation}
Finally using the master triangle $\hat{T}$ from standard arguments we know that the maximum in \eqref{boundeij} is  
bounded above by a mesh-independent constant times $h_T^{-1}$. In short we have $|e_{ij}| \leq C_E h_T$ $\forall \; i,j$, where $C_E$ is a mesh-independent constant. Hence the matrix $\tilde{K}$ equals the identity matrix plus an $O(h_T)$ matrix $E_K$. Therefore $\tilde{K}$ is an invertible matrix, as long as $h$ is sufficiently small. \QED \\ 
  
Now let us set the problem associated with spaces $V_h$ and $W_h$, whose solution is an approximation of $u$, that is, the solution of (\ref{Poisson}). Extending 
$f$ in $\Omega_h \setminus \Omega$ in different ways to be specified hereafter, and still denoting the resulting function defined in $\tilde{\Omega}_h$ by $f$, 
we wish to solve,
\begin{equation}
\label{Poissonh}
\left\{
\begin{array}{l}
\mbox{Find } u_h \in W_h \mbox{ such that } a_h(u_h,v) = F_h(v) \; \forall v \in V_h \\
\mbox{where} \\
a_h(w,v) := \int_{\Omega_h} [\nu {\bf grad}\; w \cdot {\bf grad}\; v + ({\bf b} \cdot {\bf grad}\;w) v] \mbox{ and } F_h(v) : = \int_{\Omega_h} f v.
\end{array}
\right.
\end{equation}  

We next prove: 
\begin{e-proposition}
\label{propo1}
Provided $h$ is sufficiently small problem (\ref{Poissonh}) has a unique solution. Moreover there exists a constant $\alpha > 0$ independent of $h$ such that,
\begin{equation}
\label{infsup}
\forall w \in W_h \neq 0, \displaystyle \sup_{v \in V_h \setminus \{ 0 \}} \frac{a_h(w,v)}{\parallel {\bf grad} \; w \parallel_{0,h} \parallel {\bf grad} \; v \parallel_{0,h}} 
\geq \alpha.
\end{equation}  
\end{e-proposition}

\prov Given $w \in W_h$ let $v \in V_h$ coincide with $w$ at all lattice points of elements $T \in {\mathcal T}_h \setminus {\mathcal S}_h$. As for an element $T \in {\mathcal S}_h$ we set $v=w$ at the lattice points of $T$ not belonging to $\Gamma_h$ and $v=0$ at the lattice points located on $\Gamma_h$. 
The fact that on the edges common to two elements $T^{-}$ and $T^{+}$ in ${\mathcal T}_h$, both $v_{|T^{-}}$ and $v_{|T^{+}}$ are polynomials of degree less than or equal to $k$ 
in terms of one variable coinciding at the exact number of points required to uniquely define such a function, implies indeed that $v$ is continuous in $\Omega_h$. Moreover for the same reason $v$ vanishes all over $\Gamma_h$. 
 \\
\noindent For $T \in {\mathcal S}_h$ we denote by ${\mathcal L}_T$ the set of $k-1$ lattice points of $T$ different from vertexes that belong to $\Gamma_h$. We also denote by ${\bf n}_h$ the unit outer normal vector along $\Gamma_h$. Since $div \; {\bf b} \equiv 0$ by assumption, integration by parts easily yields 
$\int_{\Omega_h} ({\bf b} \cdot {\bf grad} \; w) w = \displaystyle \oint_{\Gamma_h} \frac{{\bf b} \cdot {\bf n}_h}{2} w^2 $. Hence,   
\begin{equation}
\label{ahwv}
\left\{
\begin{array}{l} 
a_h(w,v) = \displaystyle \sum_{T \in {\mathcal T}_h} \int_T \nu |{\bf grad} \; w |^2 \\
- \displaystyle \sum_{T \in {\mathcal S}_h} \left\{ \int_T \left[ 
\nu {\bf grad} \; w \cdot {\bf grad} \; r_T(w) + ({\bf b} \cdot {\bf grad} \; w ) r_T(w) \right] - \displaystyle \int_{e_T} \frac{{\bf b} \cdot {\bf n}_h}{2} w^2 \right\},\\
\mbox{where} \\
 r_T(w) = \displaystyle \sum_{M \in {\mathcal L}_T} w(M) \varphi_M,
\end{array} 
\right.
\end{equation}
\noindent $\varphi_M$ being the canonical basis function of the space ${\mathcal P}_k(T)$ associated with the lattice point $M$. \\
\noindent Now from standard results it holds for two mesh-independent constants $C_{\varphi,0}$ and $C_{\varphi,1}$:
\begin{equation}
\label{phiestim}
\left\{
\begin{array}{l}
\parallel \varphi_M \parallel_{0,T} \leq C_{\varphi,0} h_T, \\
\\
\parallel {\bf grad} \; \varphi_M \parallel_{0,T} \leq C_{\varphi,1}.
\end{array}
\right.
\end{equation} 
On the other hand, since $w(P)=0$, where $P$ is the point of $\Gamma$ corresponding to $M \in \Gamma_h$ in accordance with the definition of $W_h$, a simple Taylor expansion about $P$ allows us to conclude that $|w(M)| \leq length(\overline{PM}) \parallel {\bf grad} \; w \parallel_{0,\infty,T^{\Delta}}$. Recalling Proposition \ref{prop01} together with \eqref{LinftyTDelta} we have, $|w(M)| \leq 
{\mathcal C}_{\infty} C_{\Gamma} h_T^2 \parallel {\bf grad} \; w \parallel_{0,\infty,T}$ for every lattice point $M \in e_T$. Then using 
\eqref{L2TDelta} it follows that,
\begin{equation}
\label{inverseineq}
|w(M)| \leq {\mathcal C}_{J} {\mathcal C}_{\infty} C_{\Gamma} h_T \parallel {\bf grad} \; w \parallel_{0,T}.
\end{equation}
On the other hand for every $Q \in e_T$ we have $|w(Q)| \leq length(\overline{OQ}) \parallel {\bf grad} \; w \parallel_{0,\infty,T}$, where $O$ is an end-point of $e_T$. Therefore it holds, 
\begin{equation}
\label{wqestim}
|w(Q)| \leq  h_T \parallel {\bf grad} \; w \parallel_{0,\infty,T} \forall Q \in e_T. 
\end{equation}
Thus from \eqref{wqestim} we infer that, 
\begin{equation*}
\int_{e_T} {\bf b} \cdot {\bf n}_h w^2 \leq h_T^3 \parallel {\bf b} \parallel_{0,\infty} \parallel {\bf grad} \; w \parallel_{0,\infty,T}^2.
\end{equation*}
Now using \eqref{L2TDelta} this further yields, 
\begin{equation}
\label{ointegestim}
\int_{e_T} {\bf b} \cdot {\bf n}_h w^2 \leq {\mathcal C}_J^2 h_T \parallel {\bf b} \parallel_{0,\infty} \parallel {\bf grad} \; w \parallel_{0,T}^2.
\end{equation}
Hence noticing that $card({\mathcal L}_T) = k-1$ $\forall T$, plugging \eqref{phiestim}, \eqref{inverseineq} and \eqref{ointegestim} into (\ref{ahwv}), we derive:
\begin{equation}
\label{ahwvbound} 
\begin{array}{l}
a_h(w,v) \geq \int_{\Omega_h} \nu |{\bf grad} \; w |^2 - h 
\displaystyle \left\{ \left[ {\mathcal C}_J^2/2 + (k-1) {\mathcal C}_{\infty} {\mathcal C}_{J} C_{\Gamma} C_{\varphi,0} \right] \parallel {\bf b} \parallel_{0,\infty} \right.\\
\left. + (k-1) {\mathcal C}_{\infty} {\mathcal C}_{J} C_{\Gamma} C_{\varphi,1} \nu \right\}   
\displaystyle \sum_{T \in {\mathcal S}_h} \parallel {\bf grad} \; w \parallel_{0,T}^2. 
\end{array}
\end{equation}         
From \eqref{ahwvbound} we readily obtain for two suitable mesh-independent constants $C_0$ and $C_1$:
\begin{equation}
\label{ahwvbelow} 
a_h(w,v) \geq [ \nu(1 - C_1 h) - C_0 \parallel {\bf b} \parallel_{0,\infty} h] \parallel {\bf grad} \; w \parallel_{0,h}^2
\end{equation}
Now using arguments in all similar to those employed above, we easily infer that,  
\begin{equation}
\label{normbound}
 \parallel {\bf grad} \; v \parallel_{0,h} \leq \parallel {\bf grad} \; w  \parallel_{0,h} + \parallel {\bf grad}(v-w) \parallel_{0,h} \leq (1+ C_1 h) 
\parallel {\bf grad} \; w \parallel_{0,h}.
\end{equation}
Combining (\ref{ahwvbelow}) and (\ref{normbound}), provided $h \leq \displaystyle \min[(4C_1)^ {-1},(4 C_0 \mbox{P\' e})^{-1}]$, where 
P\'e := $\parallel {\bf b} \parallel_{0,\infty}/\nu$ is the P\'eclet number, we establish (\ref{infsup}) with $\alpha = 2 \nu / 5$. \\

\noindent Since obviously $dim(V_h) = dim(W_h)$, the simple fact that (\ref{infsup}) holds implies that (\ref{Poissonh}) is uniquely solvable (cf. \cite{COAM}). \QED \\

\section{Error estimates}

In this section we establish errors estimates for problem \eqref{Poissonh}. To begin with we supply some useful material.  

\subsection{Preliminaries}

Several developments in the sequel rely on classical inverse inequalities applying to polynomials defined in $T \in {\mathcal T}_h$ 
(see e. g. \cite{Verfuerth}) and their extensions to neighboring sets.  
Besides \eqref{L2TDelta} we shall use the following one: \\
There exists a constant ${\mathcal C}_I$ depending only on 
$\Gamma$, $k$ and the shape regularity of ${\mathcal T}_h$ (cf. \cite{BrennerScott}, Ch.4, Sect 4) such that for any 
set $\tilde{T}$ satisfying $T^{'} \subseteq \tilde{T} \subseteq T^{\Delta}$, it holds for $1 \leq j \leq k$:  
\begin{equation}
\label{inverse}
\parallel D^j w \parallel_{0,\tilde{T}} \leq {\mathcal C}_I h_T^{-1} \parallel D^{j-1} w \parallel_{0,\tilde{T}} \; 
\forall w \in {\mathcal P}_k(\tilde{T}).
\end{equation}

Now $w$ being a given function in $H^2(\tilde{\Omega}_h)$ we define its $W_h$-interpolate $I_h(w)$ in $C^0(\tilde{\Omega}_h)$ in the following fashion:\\
Since $w$ belongs to $H^2(\Omega)$ it is possible to uniquely define $w(Q)$ at any point $Q \in \bar{\Omega}$ (cf. \cite{Adams}). In every $T \in 
{\mathcal T}_h \setminus {\mathcal S}_h$, $I_h(w)$ is the standard ${\mathcal P}_k$-interpolate of $w$ at the lattice points of $T$. If $T \in {\mathcal S}_h$, $I_h(w)$ is the unique function in ${\mathcal P}_k(T)$ such that $I_h(w)$ at $Q$ equals $w(Q)$ for all $Q$ in the union of the set of $m_k+2$ lattice points of $T$ that do not lie in the interior of $e_T$ with the set of $k-1$ points in $\Gamma$ associated with the lattice points $M$ of $T$ lying in the interior of $e_T$, as described in item 4. of the definition of $W_h$. Finally $I_h(w)$ is extended to $\bar{\Omega} \setminus \bar{\Omega}_h$ in the way prescribed for any function in $W_h$.  \\

Hereafter we will also need the following technical lemmata.

\begin{lemma}
\label{interperror}
Let $m$ be an integer, $m>1$. There exists a mesh-independent constant $C_{\Omega}$ such that $\forall w \in H^{m}(\Omega)$ such that $w_{|\Gamma} = 0$, for $j=0,1,\ldots,m$ and $p \in [1,\infty]$ it holds:
\begin{equation}
\label{interpolerror}
\parallel D^j[w - I_h(w)] \parallel_{0,p} \leq C_{\Omega} h^{m-j} | w |_{m,p}.
\end{equation}
\end{lemma}

\prov From standard results (see e.g. \cite{Ciarlet}) we know that 
\begin{equation}
\label{errest2}
\parallel D^j[w - I_h(w)]_{|T} \parallel_{0,p,T} \leq C^{'}_{\Omega} h^{m-j} | w |_{m,p,T} \; \forall T \in {\mathcal T}_h \setminus {\mathcal S}_h,
\end{equation}
where $C^{'}_{\Omega}$ is a constant independent of $h$ and $w$. \\
Now if $T \in {\mathcal S}_h$ we consider the mapping ${\mathcal G}_T$ from $T^{\Gamma}$ onto a unit element $\check{T}^{\Gamma}$ of a reference plane with coordinates 
$(\check{x},\check{y})$ given by ${\mathcal G}_T(x,y) = (x,y)/h_T$. 
From \textit{Assumption}$^{+}$ $T^{\Gamma}$ is star-shaped with respect to a ball contained in $T$.
It follows that we can extend the well-known results for the Lagrange interpolation with the set of lattice points to the one constructed in accordance with the definition of $W_h$. More precisely we mean the set consisting of the $m_k+2$ transformations in $\check{T}^{\Gamma}$ under ${\mathcal G}_T$ of lattice points of $T$  
which do not lie in the interior of $e_T$, completed with the transformations under ${\mathcal G}_T$ of the $k-1$ points $P \in \Gamma \cap T^{\Gamma}$ associated with the lattice points $M$ of $T$ lying in the interior of $e_T$ (see Figure 2). Let us denote by $\check{w}$ the transformations under ${\mathcal G}_T$ in $\check{T}^{\Gamma}$ of $w$ restricted to $T^{\Gamma}$. Notice that the transformation under  ${\mathcal G}_T$ of $I_h(w)$ is a ${\mathcal P}_k$-interpolate $\check{I}(\check{w})$ of $\check{w}$ in $\check{T}^{\Gamma}$, both functions coinciding 
whenever $\check{w}$ belongs to ${\mathcal P}_k(\check{T}^{\Gamma})$. 
Therefore, denoting by $\check{D}^j \check{v}$ the $j$-th order tensor of partial derivatives of order $j$ of a function $\check{v}$ defined in $\check{T}^{\Gamma}$, there exists a constant $\check{C}_{\mathcal T}$ depending on $\check{T}^{\Gamma}$, and hence on $\Gamma$ but not on $T$, such that,
\begin{equation}
 \label{errest1bis}
\parallel \check{D}^j[\check{w} - \check{I}(\check{w})] \parallel_{0,p,\check{T}^{\Gamma}} \leq \check{C}_{\mathcal T}  | \check{w} |_{m,p,\check{T}^{\Gamma}}.
\end{equation}
Thus, denoting by $\rho_T$ the radius of the circle inscribed in $T$, by the same arguments as in Theorem 4.4.4 of \cite{BrennerScott}, we immediately conclude that 
\begin{equation}
\label{errest1}
 \parallel D^j[w - I_h(w)]_{|T^{\Gamma}} \parallel_{0,p,T^{\Gamma}} \leq C_{\mathcal T} h_T^{m-j} 
\displaystyle \frac{h_T^{2/p}}{\rho_T^{2/p}} | w |_{m,p,T^{\Gamma}} \; \forall T \in {\mathcal S}_h, 
\end{equation} 
\noindent $C_{\mathcal T}$ being a constant depending only on $m$, $j$, $p$ and  $\check{T}^{\Gamma}$. Actually $\check{T}^{\Gamma}$ varies 
with $T$, but the underlying dependence of $C_{\mathcal T}$ reduces to a dependence on $\Gamma$ rather than on $T$ itself. \\
Now recalling that the chunkiness parameter $\sigma = \max_{T \in {\mathcal T}_h} h_T/\rho_T$ (cf. \cite{BrennerScott}) is bounded above for every ${\mathcal T}_h$ in the family of partitions in use, we set $\tilde{C}_{\Omega}:= C_{\mathcal T}\sigma^{2/p}$.\\
Finally, putting together \eqref{errest2} and \eqref{errest1} we establish \eqref{interpolerror} with $C_{\Omega} = \max[\tilde{C}_{\Omega},
C^{'}_{\Omega}]$. \QED \\

Incidentally we observe that standard approximation results (cf. \cite{BrennerScott} Subection 4.4.1) allow us to write for 
$w \in \tilde{\Omega}_h$ and a constant $\bar{C}_{\Omega}$ independent of $h$:
\begin{equation}
\label{interpolerrorh}
\parallel D^j[w - I_h(w)] \parallel_{0,h} \leq \bar{C}_{\Omega} h^{m-j} \parallel D^m w \parallel_{0,h}.
\end{equation}

\begin{lemma}
\label{wh}
Let $r=1/2+\epsilon$ for a certain $\epsilon$ in $(0,1/2)$ and $w \in H^{k+1+r}(\tilde{\Omega}_h)$ be such that 
$w_{|\Gamma} \equiv 0$. Let $\tilde{T}$ be a closed set fulfilling $T^{'} \subseteq \tilde{T} \subseteq T^{\Delta}$ and $w_h$ be a function in $W_h$ extended to $\Delta_T$ $\forall T \in {\mathcal S}_h$, as prescribed 
in Section 3. Then there exist constants ${\mathcal C}_j$ independent of $T$ and $h$ such that for $j=1,2,\ldots,k$ it holds,
\begin{equation}
\label{Djwh}
 \parallel \! D^j(w_h-w) \! \parallel_{0,\infty,\tilde{T}} \leq 
{\mathcal C}_j h_T^{-j} [\parallel\! {\bf grad}(w_h-w)\! \parallel_{0,T^{'}} + h_T^k | w |_{k+1,T^{'}} + h_T^{k+r}\! \parallel w \parallel_{k+1+r,T^{\Delta}}].
\end{equation}
\end{lemma}
      
\prov First of all recalling the interpolate $I_h(w)$ of $w$ in $\tilde{\Omega}_h$   
we write $w_h-w = (w_h-I_h(w))+(I_h(w)-w)$. Now $\forall T \in {\mathcal S}_h$ we  
make use of \eqref{L2TDelta} to write,
\begin{equation}
\label{wh1}
 \parallel D^j(w_h-w) \parallel_{0,\infty,\tilde{T}} \leq {\mathcal C}_J h_T^{-1} 
 \parallel D^j(w_h-I_h(w)) \parallel_{0,T^{'}} + \parallel D^j(I_h(w)-w) \parallel_{0,\infty,T^{\Delta}}.
\end{equation}
Applying the proportionality transformation ${\mathcal G}_T$ and setting 
$\check{T}^{\Delta}:= {\mathcal G}_T(T^{\Delta})$ we observe that 
$H^{k+1+r}(\check{T}^{\Delta})$ is continuously embedded in $W^{k,\infty}(\check{T}^{\Delta})$ (cf. \cite{Adams}), that is, there exists a constant $\check{C}_e$ depending only on $\check{T}^{\Delta}$ and hence on $\Gamma$ but not on $T$, such that,
\begin{equation}
\label{hatembed}
\parallel \check{v} \parallel_{k,\infty,\check{T}^{\Delta}} \leq \check{C}_e \parallel \check{v} \parallel_{k+1+r,\check{T}^{\Delta}} \; \forall 
\check{v} \in H^{k+1+r}(\check{T}^{\Delta}).
\end{equation}
On the other hand it is easy to see that an estimate of the type \eqref{errest1bis} also applies if $\check{T}^{\Gamma}$ is replaced by $\check{T}^{\Delta}$, with an eventual adjustment of  $\check{C}_{\mathcal T}$.
Thus taking $p= \infty$ and $m=k$ we come up with,  
\begin{equation}
\label{wh1bis}
h_T^j \parallel D^j[w - I_h(w)] \parallel_{0,\infty,T^{\Delta}} = \parallel \check{D}^j[\check{w} - \check{I}(\check{w})] \parallel_{0,\infty,\check{T}^{\Delta}} \leq \check{C}_e \check{C}_{\mathcal T} 
\parallel \check{w} \parallel_{k+1+r,\check{T}^{\Delta}}
\end{equation}
Combining \eqref{wh1bis} with standard transformation results applying to functions in fractional Sobolev spaces 
(cf. \cite{Arcangeli}), we obtain for suitable mesh-independent constants $C_j^{\Delta}$,
\begin{equation}
\label{wh2}
 \parallel D^j(I_h(w)-w) \parallel_{0,\infty,T^{\Delta}} \leq C_j^{\Delta} h_T^{k-j+r} \parallel w \parallel_{k+1+r,T^{\Delta}} \; \mbox{ for }  j=1,2,\ldots,k.
\end{equation}
On the other hand using \eqref{inverse} we easily come up with,
\begin{equation}
\label{wh3}
\parallel D^j(w_h-I_h(w)) \parallel_{0,T^{'}} \leq 
[{\mathcal C}_{I} h_T]^{-j+1} [\parallel {\bf grad}(w_h-w)) \parallel_{0,T^{'}} + \parallel {\bf grad}(w-I_h(w)) \parallel_{0,T^{'}}].
\end{equation}
Moreover using \eqref{errest1} with $p=2$ and $m=k+1$, together with standard approximation results in $T$ (cf. \cite{BrennerScott}), akin to the proof of Lemma 
\ref{interperror} we infer the existence of a mesh-independent constant ${\mathcal C}_L$ such that  
\begin{equation}
\label{wh4}
\parallel D^j(w-I_h(w)) \parallel_{0,T^{'}} \leq {\mathcal C}_L h_T^{k+1-j} | w |_{k+1,T^{'}} \mbox{ for } 1 \leq j \leq k.
\end{equation}
The combination of \eqref{wh1}, \eqref{wh2}, \eqref{wh3} and \eqref{wh4} with $j=1$ 
readily yields \eqref{Djwh} with suitable constants ${\mathcal C}_j$. \QED \\

In order to derive error estimates for problem (\ref{Poissonh}) we resort to the approximation theory of non coercive linear variational problems (cf. 
\cite{Brezzi} and \cite{COAM}). At this point it is important to recall that the solution $u$ of (\ref{Poisson}) with $g \equiv 0$ 
satisfies $a(u,v) = F(v)$ $\forall v \in H^1_0(\Omega)$, with $u \in H^1_0(\Omega)$, where, 
\begin{equation}
\label{aF}
a(w,v) := \int_{\Omega} [\nu {\bf grad} \; w \cdot {\bf grad}\; v + ({\bf b} \cdot {\bf grad} \; w ) v ] \mbox{ and } F(v) : = \int_{\Omega} fv.
\end{equation} 
Hence, owing to the construction of $V_h$, if $\Omega$ is convex $u$ also fulfills $a_h(u,v) = F_h(v) \; \forall v \in V_h$. In case $\Omega$ is not convex, 
we could extend $u$ by zero in $\Omega_h \setminus \Omega$, to define $a_h(u,v)$. However in this case there will be a non zero residual 
$a_h(u,v)-F_h(v)$ for $v \in V_h$ whose order may erode the one the approximation method (\ref{Poissonh}) is supposed to attain. Nevertheless in this case such an effect can be neutralized by means of a trick to be explained later on.

\subsection{The convex case} 

First we have,

\begin{theorem}
\label{theorem1}
As long as $h$ is sufficiently small, if $\Omega$ is convex and the solution $u$ of (\ref{Poisson}) for $g \equiv 0$ belongs to $H^{k+1}(\Omega)$, the solution $u_h$ of (\ref{Poissonh}) satisfies for $k>1$ and a suitable constant ${\mathcal C}$ independent of $h$ and $u$:
\begin{equation}
\label{errestconvex}
\parallel {\bf grad}(u - u_h) \parallel_{0,h} \leq {\mathcal C} h^k | u |_{k+1}.
\end{equation}
\end{theorem}

\prov From (\ref{infsup}) we infer that

\begin{equation}
\label{bound1} 
\parallel {\bf grad}[u_h - I_h(u)] \parallel_{0,h} \leq \displaystyle \alpha^{-1} \sup_{v \in V_h \setminus \{0\}} \frac{a_h(u_h-I_h(u),v)}{\parallel {\bf grad} \; 
v \parallel_{0,h}}.
\end{equation} 
Let us add and subtract $u$ in the first argument of $a_h$ and resort to the Friedrichs-Poincar\'e inequality, according to which $\parallel v \parallel_{0} 
\leq C_P \parallel {\bf grad}\; v \parallel_{0}$, where $C_P$ is a constant depending only on $\Omega$. In doing so we obtain after a straightforward calculation:
\begin{equation}
\label{bound2} 
\parallel {\bf grad}[u_h - I_h(u)] \parallel_{0,h} \leq \displaystyle \alpha^{-1} \left[ A \parallel {\bf grad}[u - I_h(u)] \parallel_{0,h} + 
\displaystyle \sup_{v \in V_h \setminus \{0\}} \frac{a_h(u_h-u,v)}{\parallel {\bf grad} \; v \parallel_{0,h}} \right], 
\end{equation}
where $A := \nu + C_P \parallel {\bf b} \parallel_{0,\infty}$.
Noting that $a_h(u_h,v)=F_h(v)$ we come up with:
\begin{equation}
\label{bound3} 
\parallel {\bf grad}[u_h - I_h(u)] \parallel_{0,h} \leq \displaystyle \frac{1}{\alpha} \left\{ A \parallel {\bf grad}[u - I_h(u)] \parallel_{0,h} + 
\displaystyle \sup_{v \in V_h \setminus \{0\}} \frac{|a_h(u,v)-F_h(v)|}{\parallel {\bf grad}\; v \parallel_{0,h}} \right\}.
\end{equation}
Since $\Omega_h \subset \Omega$ if $\Omega$ is convex, we observe that $a_h(u,v) = \displaystyle \oint_{\Gamma_h} \nu v \displaystyle \frac{\partial u}{\partial n_h} +  \int_{\Omega_h} v (-\nu \Delta u + {\bf b} \cdot {\bf grad} \; u)$, where $\displaystyle \frac{\partial u}{\partial n_h}$ is the outer normal derivative of $u$ on $\Gamma_h$. From equation \eqref{Poisson} and since $v \equiv 0$ on $\Gamma_h$, it trivially follows that,
\begin{equation}
\label{bound4} 
\parallel {\bf grad}(u_h - u) \parallel_{0,h} \leq \displaystyle \left(1+ \frac{A}{\alpha} \right) \parallel {\bf grad}[u - I_h(u)] \parallel_{0,h}.
\end{equation}
Finally combining (\ref{bound4}) and (\ref{interpolerror}) with $j=1$, $m=k+1$ and $p=2$, we establish \eqref{errestconvex} with 
${\mathcal C}:= [1+A/\alpha]C_{\Omega}$. \QED \\

\indent $O(h^{k+1})$-error estimates in the $L^2$-norm can be established in connection with 
Theorem \ref{theorem1}, if we require a little more regularity from $u$, according to,
\begin{theorem}
\label{theorem1bis}
As long as $h$ is sufficiently small, if $\Omega$ is convex and the solution $u$ of (\ref{Poisson}) for $g \equiv 0$ belongs to $H^{k+1+r}(\Omega)$ with $r=1/2+ \epsilon$ for $\epsilon > 0$ arbitrarily small, the solution $u_h$ of (\ref{Poissonh}) satisfies for $k>1$ and a suitable constant ${\mathcal C}_0$ independent of $h$ and $u$:
\begin{equation}
\label{L2estconvex}
\parallel u - u_h \parallel_{0,h} \leq {\mathcal C}_0 h^{k+1} \parallel u \parallel_{k+1+r}.
\end{equation}
\end{theorem}

\prov Recalling that every function in $W_h$ is defined in $\bar{\Omega} \setminus \Omega_h$, let $\bar{u}_h$ be the function given by $\bar{u}_h = u_h-u$ in $\Omega$. Let also $v \in H^1_0(\Omega)$ be the solution of 
\begin{equation}
\label{adjoint}
- \nu \Delta v - {\bf b} \cdot {\bf grad}\; v = \bar{u}_h \; \mbox{in } \Omega.
\end{equation}
Since both $\Omega$ and ${\bf b}$ are sufficiently smooth and $\bar{u}_h \in L^2(\Omega)$ we know that $v \in H^2(\Omega)$ and moreover there 
exists a constant $C(\Omega)$ depending only on $\nu$, ${\bf b}$ and $\Omega$ such that,
\begin{equation}
\label{adjoint1}
\parallel v \parallel_{2} \leq C(\Omega) \parallel \bar{u}_h \parallel_{0}. 
\end{equation}
Therefore
\begin{equation}
\label{L2est1} 
\parallel \bar{u}_h \parallel_{0} \leq C(\Omega) \displaystyle \frac{\int_{\Omega} \bar{u}_h (- \nu \Delta v - {\bf b} \cdot {\bf grad}\; v)}{\parallel v \parallel_{2}}.
\end{equation}
Using integration by parts we easily obtain,
\begin{equation}
\label{L2est2} 
\parallel \bar{u}_h \parallel_{0} \leq C(\Omega) \displaystyle \frac{a(\bar{u}_h, v)
+ b_{1h}(\bar{u}_h,v)}{\parallel v \parallel_{2}} 
\end{equation}
where
\begin{equation}
\label{b1h}
b_{1h}(w,v):= -\nu \displaystyle \int_{\Gamma} w  \frac{\partial v}{\partial n} \mbox{ for } w \in H^1(\Omega) 
\mbox{ and } v \in H^1_0(\Omega) \cap H^2(\Omega).
\end{equation}
Let $\Pi_h(v)$ be the continuous piecewise linear interpolate of $v$ in $\Omega$ at the vertexes of the mesh. 
Setting $v_h=\Pi_h(v)$ in $\Omega_h$ we have $v_h \in V_h$. Therefore it holds $a(u,v_h)=a_h(u,v_h)=F(v_h)=F_h(v_h)=a_h(u_h,v_h)$. On the other hand 
$a(\bar{u}_h,v)=a_h(\bar{u}_h,v) + a_{\Delta_h}(\bar{u}_h,v)$ where 
\begin{equation}
\label{aDeltah}
a_{\Delta_h}(w,z):=\int_{\Delta_h} [\nu {\bf grad}\; w \cdot {\bf grad}\;z + {\bf b} \cdot {\bf grad}\; w \; z] \mbox{ for } w,z \in H^1(\Omega) \mbox{ with } \Delta_h = \Omega \setminus \Omega_h.
\end{equation}
Now we observe that $a_{\Delta_h}(\bar{u}_h,v)=a_{\Delta_h}(\bar{u}_h,v-\Pi_h(v))
+a_{\Delta_h}(\bar{u}_h,\Pi_h(v))$. Thus applying First Green's identity in $\Delta_T$ for $T \in {\mathcal S}_h$we come up with,
$a_{\Delta_h}(\bar{u}_h,\Pi_h(v)) = b_{2h}(\bar{u}_h,\Pi_h(v))+b_{3h}(\bar{u}_h,\Pi_h(v))$, where
\begin{equation}
\label{b2h}
b_{2h}(w,z):=\displaystyle \sum_{T \in {\mathcal S}_h} \int_{\Delta_T} [-\nu \Delta w+{\bf b} \cdot {\bf grad} \; w] z \mbox{ for } w \in W_h + H^2(\Omega) \mbox{ and } z \in H^{1}(\Omega),  
\end{equation}
and setting $\partial T = T^{\Delta} \cap \Gamma$ for $T \in {\mathcal S}_h$,  
\begin{equation}
\label{b3h}
b_{3h}(w,z) := \nu \displaystyle \sum_{T \in {\mathcal S}_h} \int_{\partial T} \frac{\partial w}{\partial n} z \mbox{ for } w \in W_h + H^2(\Omega) \mbox{ and } z \in H^1(\Omega).
\end{equation}
Further setting
\begin{equation}
\label{b4h} 
b_{4h}(w,z) := a_{\Delta_h}(w,z) \mbox{ for } w,z \in H^1(\Omega),  
\end{equation}  
it follows that,
\begin{equation}
\label{L2est3} 
\begin{array}{l}
\parallel \bar{u}_h \parallel_{0} \leq C(\Omega) \displaystyle \frac{a_h(\bar{u}_h,e_h(v))
+ b_{1h}(\bar{u}_h,v)+b_{2h}(\bar{u}_h,\Pi_h(v))+b_{3h}(\bar{u}_h,\Pi_h(v))+
b_{4h}(\bar{u}_h,e_h(v))}{\parallel v \parallel_{2}}, \\
\mbox{with } e_h(v) = v -\Pi_h(v). 
\end{array}
\end{equation} 
Now from classical results, for a mesh-independent constant $C_{V}$ it holds
\begin{equation}
\label{interP1}
\parallel {\bf grad} \; e_h(v) \parallel_{0,h} \leq \parallel {\bf grad} \; e_h(v) \parallel_{0} \leq C_{V} h 
| v |_{2}.
\end{equation} 
Since $a_h(\bar{u}_h,e_h(v)) \leq A \parallel {\bf grad}\; \bar{u}_h \parallel_{0,h} \parallel {\bf grad}\; e_h(v) \parallel_{0,h}$, using \eqref{errestconvex} and \eqref{interP1} we obtain,
\begin{equation}
\label{L2est3bis}
a_h(\bar{u}_h,e_h(v) ) \leq A C_{V} {\mathcal C} h^{k+1} |u|_{k+1} |v|_2.
\end{equation}
Hence, setting $\tilde{\mathcal C}_0:= C(\Omega)A C_{V} {\mathcal C}$ and recalling \eqref{L2est3}, it holds,
\begin{equation}
\label{L2est4} 
\parallel \bar{u}_h \parallel_{0} \leq \tilde{\mathcal C}_0 h^{k+1} | u |_{k+1} + C(\Omega) \displaystyle \frac{b_{1h}(\bar{u}_h,v)+b_{2h}(\bar{u}_h,\Pi_h(v))+b_{3h}(\bar{u}_h,\Pi_h(v))+
b_{4h}(\bar{u}_h,e_h(v))}{\parallel v \parallel_{2}}.
\end{equation} 

Let us estimate $b_{ih}$ for $i=1,2,3,4$.\\
As for $b_{1h}$ we first note that according to the Trace Theorem there exists a constant $C_t$ depending only on $\Omega$ such that 
\begin{equation}
\label{L2estimate}
b_{1h}(\bar{u}_h,v) \leq C_t \parallel \bar{u}_h \parallel_{0,\Gamma} \parallel v \parallel_{2}.
\end{equation} 
Now for every $T \in {\mathcal S}_h$ we take a local orthogonal frame $(O;x,y)$ whose origin $O$ is a vertex of $T$ in $\Gamma$, $x$ is the abscissa along the edge $e_T$. Let $s$ be the curvilinear abscissa along $\partial T$ with origin at $O$. Notice that owing to our assumptions $s$ can be uniquely expressed in terms of $x$ and conversely, for $x \in [0,l_T]$, where $l_T$ is the length of $e_T$. Adopting $f_T(x)$ as the $y$-abscissa of the points in $\partial T$, let $\breve{u}_h$ be the function of $x$ defined by $\breve{u}_h(x) = \bar{u}_h[x,f_T(x)]$. Since $\breve{u}_h$ vanishes at $k+1$ different points in $[0,l_T]$, from standard results for one-dimensional interpolation (cf \cite{Quarteroni}), there exists a mesh-independent constant $C_L$ such that,
\begin{equation}
\label{auxiliary0}
\displaystyle \left[\int_0^{l_T} |\breve{u}_h(x)|^2 dx \right]^{1/2} \leq C_L h_T^{k+1} \displaystyle \left[\int_0^{l_T} \left|[\breve{u}_h^{(k+1)}](x) \right|^2 dx \right]^{1/2}. 
\end{equation}
On the other hand, recalling Proposition \ref{prop02}, there exist mesh-independent constants $c_{j,\Gamma}$ such that,
\begin{equation}
\label{auxiliary1}
\displaystyle \max_{x \in [0,l_T]}  |f_T^{(j)}(x)| \leq c_{j,\Gamma} h_T^{2-j}, \; j=1,2,\ldots,k+1 \; \forall T \in {\mathcal S}_h.
\end{equation}
Thus taking into account that the derivatives of $u_h$ of order greater than $k$ vanish in $T^{\Delta}$, straightforward calculations using the chain rule yield for suitable mesh-independent constants $c_{j}$, $j=0,1,\ldots,k$: \\
\begin{equation}
\label{auxiliary2}
|\breve{u}_h^{(k+1)}| \leq  c_{0} | D^{k+1}(u) | + \displaystyle 
\sum_{j=1}^{k} c_{j} h_T^{1-j} | D^{k+1-j}(\bar{u}_h) |.  
\end{equation}
All the partial derivatives appearing in \eqref{auxiliary2} are to be understood at a (variable) point in $\partial T$. \\
Now since $ds = \sqrt{1+(f_T^{'})^2}dx$ we have $length(\partial T) \leq C_q l_T$ with  $C_q:=\sqrt{1+(h_0 {\mathcal G}_{max})^2}$. Thus it holds
\begin{equation}
\label{auxiliary3}
\parallel \bar{u}_h \parallel_{0,\Gamma}^2 \leq 
C_q \displaystyle \sum_{T \in {\bf S}_h} \int_0^{l_T} |\breve{u}_h(x)|^2 dx.
\end{equation}
Combining \eqref{auxiliary0}, \eqref{auxiliary1}, \eqref{auxiliary2} and \eqref{auxiliary3}, after straightforward calculations we come up with a mesh-independent constant $C_{k,0}$ such that,
\begin{equation}
\label{auxiliary4}
\parallel \bar{u}_h \parallel_{0,\Gamma}^2 
\leq C_{k,0}  \displaystyle \left[ h^{2(k+1)} \int_{\Gamma}  |D^{k+1}(u)|^2 ds  + \sum_{T \in {\bf S}_h} h_T^{2(k+1)} 
\displaystyle \int_{\partial T} \displaystyle \sum_{j=1}^k h_T^{2(1-j)}|D^{k+1-j}(\bar{u}_h)|^2 ds  \right].\\
\end{equation}
From the Trace Theorem \cite{Adams} we know that there exists a constant $C_{r}$ such that,
\begin{equation}
\label{auxiliary5}
\int_{\Gamma}  |D^{k+1}(u)|^2  \leq 
C_{r}^2 \parallel u \parallel_{k+1+r}^2
\end{equation}
On the other hand, using again the curved triangle $T^{\Delta}$, by standard calculations we can write:
\begin{equation}
\label{auxiliary6}
\displaystyle \int_{\partial T}  
\displaystyle \sum_{j=1}^k h_T^{2(1-j)}|D^{k+1-j}(\bar{u}_h)|^2 ds 
\leq C_{q} h_T \displaystyle \sum_{j=1}^k h_T^{2(j-k)} 
\parallel D^{j}(\bar{u}_h) \parallel_{0,\infty,T^{\Delta}}^2 ds.
\end{equation} 
Now using Lemma \ref{wh} with $w=u$ and $w_h=u_h$, after straightforward calculations we obtain for a suitable mesh-independent constant 
$C_{k,1}$:
\begin{equation}
\label{auxiliary7}
\begin{array}{l}
\displaystyle \int_{\partial T} \! \displaystyle \sum_{j=1}^k \! h_T^{2-2j}|D^{k+1-j}(\bar{u}_h)|^2 \! ds  \leq 
 C_q C_{k,1} 
 h_T^{1-2k} [ \parallel\! {\bf grad}\; \bar{u}_h \! \parallel_{0,T}^2\! + h_T^{2k} | u |_{k+1,T}^2 \!+ \! h_T^{2k+2r}\! \parallel u \parallel_{k+1+r,T^{\Delta}}^2 \! ].
\end{array}
\end{equation}
From \eqref{auxiliary7} and \eqref{errestconvex} we infer that for another mesh-independent constant $C_{k,2}$ it holds: 
\begin{equation}
\label{auxiliary8}
\sum_{T \in {\bf S}_h} h_T^{2(k+1)} 
\displaystyle \int_{\partial T} \displaystyle \sum_{j=1}^k h_T^{2-2j}|D^{k+1-j}(\bar{u}_h)|^2 ds
\leq \displaystyle  C_{k,2} \left[ h^{2k+3} | u |_{k+1}^2 + h^{2k+4+2\epsilon}\! \parallel u \parallel_{k+1+r}^2 \right]. 
\end{equation}
Taking into account \eqref{auxiliary4}, \eqref{auxiliary5} and \eqref{auxiliary8}, and since $h <1$, we easily obtain,
\begin{equation}
\label{auxiliary9}
\parallel \bar{u}_h \parallel_{0,\Gamma} 
\leq C_{k,3} h^{k+1} \displaystyle \left[ | u |_{k+1} + \parallel u \parallel_{k+1+r} \right],
\end{equation} 
where $C_{k,3}$ is another mesh-independent constant.\\
It follows from \eqref{L2estimate} and \eqref{auxiliary9} that for $C_{b1} = 2C_{k,3} C_t $ it holds:
\begin{equation}
\label{estimateb1}
b_{1h}(\bar{u}_h,v) \leq C_{b1} h^{k+1} \parallel u \parallel_{k+1+r} \parallel v \parallel_{2}.
\end{equation} 

Now we turn our attention to $b_{2h}$.\\
First observing that ${\bf grad}\;\Pi_h(v)$ is constant in $T^{\Delta}$ for $T \in {\mathcal S}_h$ and $\Pi_h(v) =0$ on $\Gamma_h$, by Rolle's Theorem
\begin{equation}
\label{estim1b2}
|\Pi_h(v)(P)|\leq C_{\Gamma} h_T^2 | [{\bf grad} \; \Pi_h(v)]_{|T} | \; \forall P \in \Delta_T \mbox{ and } \forall T \in {\mathcal S}_h.
\end{equation}
Noticing that $area(\Delta_T) \leq C_{\Gamma} h_T^3$, using \eqref{estim1b2} we have,
\begin{equation}
\label{estim2b2}
b_{2h}(\bar{u}_h,\Pi_h{v}) \leq C_{\Gamma}^2 \displaystyle \sum_{T \in {\mathcal S}_h}  h_T^5 
\parallel -\nu \Delta \bar{u}_h + {\bf b} \cdot {\bf grad} \; \bar{u}_h \parallel_{0,\infty,T^{\Delta}} 
\parallel {\bf grad} \; \Pi_h(v) \parallel_{0,\infty,T}.
\end{equation}
From the classical inverse inequality (cf. \cite{Ciarlet}),
\begin{equation} 
\label{invinity}
\parallel w \parallel_{0,\infty,T} \leq C_{\infty}(k) h_T^{-1} \parallel w \parallel_{0,T} \; \forall w \in {\mathcal P}_k(T),
\end{equation}
with mesh-independent constants $C_{\infty}(k)$, we further obtain:
\begin{equation}
\label{estim3b2}
b_{2h}(\bar{u}_h,\Pi_h{v}) \leq C_{\Gamma}^2 C_{\infty}(1) \displaystyle \sum_{T \in {\mathcal S}_h}  h_T^4 
\parallel -\nu \Delta \bar{u}_h + {\bf b} \cdot {\bf grad} \; \bar{u}_h \parallel_{0,\infty,T^{\Delta}} 
\parallel {\bf grad} \; \Pi_h(v) \parallel_{0,T}.
\end{equation}
Next using Young's inequality we rewrite \eqref{estim3b2} as, 
\begin{equation}
\label{estim4b2} 
\begin{array}{l}
b_{2h}(\bar{u}_h,\Pi_h{v}) \leq  C_{\Gamma}^2 C_{\infty}(1) \displaystyle \sum_{T \in {\mathcal S}_h}  h_T^4 
[ \nu \parallel \Delta \bar{u}_h \parallel_{0,\infty,T^{\Delta}} + \parallel {\bf b} \parallel_{0,\infty} \parallel {\bf grad}\; \bar{u}_h \parallel_{0,\infty,T^{\Delta}}] \parallel {\bf grad} \; \Pi_h(v) \parallel_{0,T}.
\end{array}
\end{equation}
Thus applying Lemma \ref{wh} to \eqref{estim4b2}, after straightforward manipulations, we can assert that there exists a mesh-independent constant $\tilde{C}_2$ such that ,
\begin{equation}
\label{estim4-5b2} 
\begin{array}{l}
b_{2h}(\bar{u}_h,\Pi_h{v}) \leq \tilde{C}_2 \displaystyle \sum_{T \in {\mathcal S}_h}  h_T^2 
[ \parallel {\bf grad}\; \bar{u}_h \parallel_{0,T} + h^k  |u|_{k+1,T} + h^{k+r} \parallel u \parallel_{k+1+r,T^{\Delta}} ] \parallel {\bf grad} \; \Pi_h(v) \parallel_{0,T}.
\end{array}
\end{equation}
Applying the Cauchy-Schwarz inequality to \eqref{estim4-5b2} and recalling \eqref{errestconvex}, we obtain for another mesh-independent 
constant $\bar{C}_2$: 
\begin{equation}
\label{estim5b2} 
\begin{array}{l}
b_{2h}(\bar{u}_h,\Pi_h{v}) \leq \bar{C}_2 h^2 
[  h^k  |u|_{k+1} + h^{k+r} \parallel u \parallel_{k+1+r} ] \parallel {\bf grad} \; \Pi_h(v) \parallel_{0,h}.
\end{array}
\end{equation}
On the other hand, using \eqref{interP1} and setting $\bar{C}_V:=\sqrt{2+2h_0^2 C_{V}^2}$, we easily infer that,  
\begin{equation}
\label{estim6b2}
\parallel {\bf grad} \; \Pi_h(v) \parallel_{0,h} \leq \bar{C}_V \parallel v \parallel_{2}.
\end{equation}
Combining \eqref{estim5b2} and \eqref{estim6b2} we derive for $C_{b2}=2 \bar{C}_2 \bar{C}_V$:
\begin{equation}
\label{estimateb2}
b_{2h}(\bar{u}_h,\Pi_h{v}) \leq C_{b2} h^{k+2} \parallel u \parallel_{k+1+r} \parallel v \parallel_{2}.
\end{equation} 

Next we estimate $b_{3h}$.\\
Recalling \eqref{b3h} and the fact $\parallel {\bf grad}\;\Pi_h(v) \parallel_{0,\infty,T^{\Delta}} = \parallel {\bf grad}\;\Pi_h(v) \parallel_{0,\infty,T}$, we first define the function $\omega_T:=|[{\bf grad}\;\bar{u}_{h}]_{|T}|$ 
for every $T \in {\mathcal S}_h$. Then we have:
\begin{equation}
\label{estim1b3}
b_{3h}(\bar{u}_h,\Pi_h(v)) \leq \nu \displaystyle \sum_{T \in {\mathcal S}_h} 
\int_{\partial T} \omega_T |\Pi_h(v)| \leq \nu C_{\Gamma} h_T^2 \parallel {\bf grad} \; \Pi_h(v) \parallel_{0,\infty,T} \int_{\partial T} \omega_T.
\end{equation} 
Resorting again to the standard master triangle $\hat{T}$ we denote by $\partial\hat{T}$ the transformation of $\partial T$ 
under the affine mapping ${\mathcal F}_T$ from $T$ onto $\hat{T}$. \\ 
Recalling \eqref{estim1b2} and taking into account that $length(\partial T) \leq C_q h_T$ and $length(\partial \hat{T})=1$, we have:
\begin{equation}
\label{estim2b3}
b_{3h}(\bar{u}_h,\Pi_h(v)) \leq \nu C_{\Gamma} C_q \displaystyle \sum_{T \in {\mathcal S}_h} h_T^3 
\parallel {\bf grad} \; \Pi_h(v) \parallel_{0,\infty,T} \int_{\partial \hat{T}} \hat{\omega}, 
\end{equation}  
where $\hat{\omega}$ is the transformation of $\omega_T$ under the mapping 
${\mathcal F}_T$. \\
Next we apply the Trace Theorem to the transformation $\hat{T}^{\Delta}$ of $T^{\Delta}$ under ${\mathcal F}_T$. Thanks to the fact that $\Gamma$ is smooth and $h$ is sufficiently small, there exists a constant $\hat{C}_t$ independent of $T$ such that,
\begin{equation}
\label{estim3b3} 
\displaystyle \int_{\partial \hat{T}} \hat{\omega} \leq \hat{C}_t \displaystyle \left\{ \int_{\hat{T}^{\Delta}} [ \hat{\omega}^2 + |\widehat{\bf grad}\; \hat{\omega}|^2 ]  \right\}^{1/2},   
\end{equation}  
where $\widehat{\bf grad}$ is the gradient operator for functions defined in $\hat{T}^{\Delta}$.\\
Moving back to $T^{\Delta}$ and using \eqref{L2TDelta}, from \eqref{estim2b3} and \eqref{estim3b3} we conclude that for a suitable mesh-independent constant $\tilde{C}_{3}$ it holds,
\begin{equation}
\label{estim4b3}
b_{3h}(\bar{u}_h,\Pi_h(v)) \leq \tilde{C}_{3} \displaystyle \sum_{T \in {\mathcal S}_h} 
h_T \parallel {\bf grad} \; \Pi_h(v) \parallel_{0,T} \left\{ \int_{T^{\Delta}} [ \omega_T^2 + h_T^2|{\bf grad}\;\omega_T|^2 ] \right\}^{1/2}. 
\end{equation}
By the Cauchy-Schwarz inequality this further yields,  
\begin{equation}
\label{estim5b3}
b_{3h}(\bar{u}_h,\Pi_h(v)) \leq \tilde{C}_{3} h \parallel {\bf grad} \; \Pi_h(v) \parallel_{0,h}   
\displaystyle \left[ \sum_{T \in {\mathcal S}_h} \left( \parallel {\bf grad}\;\bar{u}_{h} \parallel_{0,T^{\Delta}}^2  + h_T^2 \parallel H(\bar{u}_{h}) \parallel_{0,T^{\Delta}}^2 \right) \right]^{1/2}. 
\end{equation}
Now we note that  
\begin{equation}
\label{estim0b3}
\parallel H(\bar{u}_h) \parallel_{0,T^{\Delta}} \leq \sqrt{area(T^{\Delta})} \parallel H(\bar{u}_h) \parallel_{0,\infty,T^{\Delta}}.
\end{equation}
Observing that  $area(T^{\Delta}) \leq area(T) + C_{\Gamma} h_T^3 \leq h_T^2 (1/2 + C_{\Gamma} h_0)$, resorting again to Lemma \ref{wh} 
we easily conclude from \eqref{estim5b3} and \eqref{estim0b3} that there is a mesh-independent constant  $\bar{C}_3$ such that 
\begin{equation}
\label{estim6b3}
\displaystyle  \sum_{T \in {\mathcal S}_h} \left( \parallel {\bf grad}\;\bar{u}_{h} \parallel_{0,T^{\Delta}}^2  + h_T^2 \parallel H(\bar{u}_{h}) \parallel_{0,T^{\Delta}}^2 \right)  \leq \bar{C}_3^2 h^{2k} [| u |_{k+1}^2 + h^{2r} \parallel u \parallel_{k+1+r}^2 ].
\end{equation}
Then plugging \eqref{estim6b2} and \eqref{estim6b3} into \eqref{estim5b3} 
and setting $C_{b3}=2 \bar{C}_3 \tilde{C}_{3} \bar{C}_V$ we come up with:
\begin{equation}
\label{estimateb3} 
b_{3h}(\bar{u}_h,\Pi_h(v)) \leq C_{b3} h^{k+1} \parallel u \parallel_{k+1,r} \parallel v \parallel_{2}. 
\end{equation}
Finally we estimate $b_{4h}$.\\
Using a geometric argument already exploited above, together with \eqref{Djwh} with $j=1$,  we can successively write:
\begin{equation}
\label{estim1b4}
b_{4h}(\bar{u}_h,v-\Pi_h(v)) \leq A C_{\Gamma}^{1/2} \displaystyle \sum_{T \in {\mathcal S}_h} h_T^{3/2} 
\parallel {\bf grad} \; \bar{u}_h \parallel_{0,\infty,T^{\Delta}} \parallel {\bf grad}( v - \Pi_h(v)) 
\parallel_{0,T^{\Delta}}, 
\end{equation}

\begin{equation}
\label{estim2b4}
\begin{array}{l}
b_{4h}(\bar{u}_h,v-\Pi_h(v)) \leq A C_{\Gamma}^{1/2} {\mathcal C}_1 \displaystyle \sum_{T \in {\mathcal S}_h} h_T^{1/2} 
[ \parallel {\bf grad} \; \bar{u}_h \parallel_{0,T} + h_T^k | u |_{k+1,T} + h_T^{k+r} \parallel u \parallel_{k+1+r,T}] \\
\parallel {\bf grad}( v - \Pi_h(v)) \parallel_{0,T^{\Delta}}, 
\end{array}
\end{equation}
Now from the Cauchy-Schwarz inequality and \eqref{errestconvex} we obtain for a mesh-independent constant $\tilde{C}_4$:
\begin{equation}
\label{estim3b4}
b_{4h}(\bar{u}_h,v-\Pi_h(v)) \leq \tilde{C}_4 h^{k+1/2} [| u |_{k+1} + h^{r} \parallel u \parallel_{k+1+r}] 
\parallel {\bf grad}( v - \Pi_h(v)) \parallel_{0}. 
\end{equation}
On the other hand using \eqref{interP1} we obtain from \eqref{estim3b4}:
\begin{equation}
\label{estim4b4}
b_{4h}(\bar{u}_h,v-\Pi_h(v)) \leq C_V \tilde{C}_4 h^{k+3/2} [| u |_{k+1} + h^{r} \parallel u \parallel_{k+1+r}]  | v |_2 
\end{equation}
Thus setting $C_{b4} = 2 \tilde{C}_4 C_{V}$ we obtain,
\begin{equation}
\label{estimateb4}
b_{4h}(\bar{u}_h,v-\Pi_h(v)) \leq C_{b4} h^{k+3/2} \parallel  u \parallel_{k+1+r} \parallel v \parallel_{2}. 
\end{equation}
 
Finally, plugging \eqref{estimateb1}, \eqref{estimateb2}, \eqref{estimateb3} and \eqref{estimateb4} into  
\eqref{L2est4}, owing to the fact that $h <1$, we immediately obtain \eqref{L2estconvex} with 
${\mathcal C}_0=\tilde{\mathcal C}_0 + C(\Omega) (C_{b1}+C_{b2}+C_{b3}+C_{b4})$. \QED \\

\subsection{The non convex case}

We next address the case of a non convex $\Omega$. Let us consider a smooth domain $\tilde{\Omega}$ close to $\Omega$ which strictly contains $\tilde{\Omega}_h$ for all $h$ sufficiently small, say $h \leq h_0$. More precisely, denoting by $\tilde{\Gamma}$ the boundary of $\tilde{\Omega}$ we assume that $meas(\tilde{\Gamma})-meas(\Gamma) \leq \varepsilon$ for $\varepsilon$ sufficiently small. Henceforth we also consider that $f$ is extended to  
$\tilde{\Omega} \setminus \Omega$. We still denote the extended function by $f$, which is arbitrarily chosen in $\tilde{\Omega} \setminus \Omega$, except for the requirement that $f \in H^{k-1}(\tilde{\Omega})$. 

Then the following theorem holds:

\begin{theorem}
\label{theorem2}
Assume that there exists a function $\tilde{u}$ defined in 
$\tilde{\Omega}$ having the properties:
\begin{itemize}
\item 
$- \nu \Delta \tilde{u} + {\bf b} \cdot {\bf grad} \; \tilde u = f$ in $\tilde{\Omega}$;
\item
$\tilde{u}_{|\Omega} = u$;
\item
$\tilde{u} = 0$ a.e. on $\Gamma$; 
\item
$\tilde{u} \in H^{k+1}(\tilde{\Omega})$.
\end{itemize}
Then as long as $h$ is sufficiently small there exists a mesh-independent constant $\tilde{\mathcal C}$ such that: 
\begin{equation}
\label{errestconcave}
\parallel {\bf grad}(u_h - u) \parallel_{0,h}^{'} \leq \tilde{\mathcal C} h^k | \tilde{u} |_{k+1,\tilde{\Omega}},
\end{equation} 
\end{theorem}

\prov First of all we extend every $v \in V_h$ by zero in $\tilde{\Omega} \setminus \Omega_h$. Then we define 
$\tilde{A}:=\nu + \tilde{C}_P \parallel b \parallel_{0,\infty,\tilde{\Omega}}$; $\tilde{C}_P$ being a constant for which the Friedrichs-Poincar\'e inequality holds in $\tilde{\Omega}$, that is, $\parallel v \parallel_{0,\tilde{\Omega}} \leq \tilde{C}_P \parallel {\bf grad} 
\; v \parallel_{0,\tilde{\Omega}} \; \forall v \in H^1_0(\tilde{\Omega})$.\\
In doing so, thanks to \eqref{infsup}, it is easy to see that an inequality analogous to \eqref{bound4} holds for $\tilde{u}$ instead of $u$, namely,
\begin{equation}
\label{boundprime}
\parallel {\bf grad} (u_h - \tilde{u}) \parallel_{0,h} \leq \displaystyle \left( 1 + \frac{\tilde{A}}{\alpha} \right) 
\parallel {\bf grad}[\tilde{u}-I_h(\tilde{u})] \parallel_{0,h}.
\end{equation}
 Since $\parallel {\bf grad}(u - u_h) \parallel_{0,h}^{'} \;
\leq \; \parallel {\bf grad}(\tilde{u} - u_h) \parallel_{0,h}$, \eqref{interpolerrorh} with $m=k+1$ and $p=2$ immediately yields (\ref{errestconcave}). \QED \\

It is noteworthy that the knowledge of a regular extension of the right hand side $f$ associated with a regular extension $\tilde{u}$ of $u$ is necessary to 
optimally solve problem (\ref{Poissonh}) in the general case. Of course, except for very particular situations such as the toy problems used to illustrate the performance of our method in \cite{arXiv2D} (see also the Appendix hereafter), in most cases such an extension of $f$ is not known. Even if we go the other around by prescribing a regular $f$ in $\tilde{\Omega}$, the existence of an associated $\tilde{u}$ fulfilling the assumptions of Theorem \ref{theorem2} can also be questioned. However using the results in \cite{Coffman} combined with standard ones (cf. \cite{LionsMagenes}) it is possible to identify cases where such an extension $\tilde{u}$ does exist. We refer to \cite{arXiv2D} for further details.\\ 
In the general case however, a convenient way to bypass the uncertain existence of an extension $\tilde{u}$ satisfying the assumptions of Theorem \ref{theorem2}, is to resort to numerical integration on the right hand side. Under certain conditions rather easily satisfied, this leads to the definition of an alternative approximate problem, in which only values of $f$ in $\Omega$ come into play. This trick is inspired by the celebrated work due to Ciarlet and Raviart on the isoparametric finite element method (cf. \cite{CiarletRaviart} and \cite{Ciarlet}). To be more specific, these authors employ the following argument, assuming that $h$ is small enough: if a numerical integration formula is used, which has no integration points different from vertexes on the edges of a triangle, then only values of $f$ in $\Omega$ will be needed to compute the corresponding approximation of $F_h(v)$. This means that the knowledge of $\tilde{u}$, and thus of the regular extension of $f$, will not be necessary for implementation purposes. Moreover, provided the accuracy of the numerical integration formula is compatible with method's order, the resulting modification of (\ref{Poissonh}) will be a method of order $k$ in the norm  $\parallel \cdot \parallel_{0,h}^{'}$ of ${\bf grad}(u - u_h)$. \\
Nevertheless it is possible to get rid of the above argument based on numerical integration in the most important cases in practice, namely, those of quadratic and cubic Lagrange finite elements. Let us see how this works. \\
First of all we consider that $f$ is extended by zero in $\Delta_{\Omega}:=\tilde{\Omega} \setminus \bar{\Omega}$, and resort to the extension $\tilde{u}$ of $u$ to the same set constructed in accordance to Stein et al. \cite{Stein}. This extension does not satisfy $- \nu \Delta \tilde{u} + {\bf b} \cdot {\bf grad} \; \tilde{u}=0$ in $\Delta_{\Omega}$ but the function denoted in the same way such that $\tilde{u}_{|\Omega} = u$ does belong to $H^{k+1}(\tilde{\Omega})$. Since $k > 1$ this means in particular that the traces of the functions $u$ and $\tilde{u}$ coincide on $\Gamma$ and that $\partial u /\partial n = - \partial \tilde{u}/\partial \tilde{n} = 0$ a.e. on $\Gamma$ where the normal derivatives on the left and right hand side of this relation are outer normal derivatives with respect to $\Omega$ and $\Delta_{\Omega}$ respectively (the trace of the Laplacian of both functions also coincide  on $\Gamma$ but this is not relevant for our purposes). Based on this extension of $u$ to $\Omega_h$ for all such polygons of interest, we next prove the following results for the approximate problem (\ref{Poissonh}), without assuming that $\Omega$ is convex. Here $f$ represents the function identical to the right hand side datum of (\ref{Poisson}) in $\Omega$, that vanishes identically in $\Delta_{\Omega}$.

\begin{theorem} 
\label{P2}
Let $k=2$ and assume that $u \in H^3(\Omega)$. Provided $h$ is sufficiently small, there exists a mesh-independent constant $C_2$ such that the unique solution $u_h$ to (\ref{Poissonh}) satisfies:
\begin{equation}
\label{estimateP2} 
\parallel {\bf grad}(u - u_h) \parallel_{0,h}^{'} \leq C_2 h^{2} G(\tilde{u})  
\end{equation}
where $\tilde{u} \in H^3(\tilde{\Omega})$ is the regular extension of $u$ to $\tilde{\Omega}$ constructed in accordance to Stein et al. \cite{Stein} and $G(\tilde{u}):=| \tilde{u} |_{3,\tilde{\Omega}} + h^{1/2} \parallel {\mathcal M}\tilde{u} \parallel_{0,\tilde{\Omega}}$, where 
${\mathcal M} \tilde{u}:= \nu \Delta \tilde{u} - {\bf b} \cdot {\bf grad} \; \tilde{u}$.
\end{theorem}

\prov First of all, straightforward manipulations lead to, 
\begin{equation}
\label{fourthbound}
\parallel {\bf grad}(u_h - w) \parallel_{0,h} \leq \displaystyle \frac{1}{\alpha} 
\left[\tilde{A}\parallel {\bf grad}(\tilde{u} - w) \parallel_{0,h} + \displaystyle \sup_{v \in V_h \setminus \{0\}} 
\frac{|a_h(\tilde{u},v)-F_h(v)|}{\parallel {\bf grad}\; v \parallel_{0,h}} \right] \; \forall w \in W_h.
\end{equation}
First of all we take $w=I_h(\tilde{u})$ and recall \eqref{interpolerrorh} to obtain,
\begin{equation}
\label{interP2}
\parallel {\bf grad}(\tilde{u} - w) \parallel_{0,h} \leq \bar{C}_{\Omega} h^2 |\tilde{u}|_{3,\tilde{\Omega}}.
\end{equation}
On the other hand the numerator in (\ref{fourthbound}) is estimated as follows: Since $\tilde{u} \in H^3(\tilde{\Omega})$ we can apply First Green's formula to $a_h(\tilde{u},v)$ thereby getting rid of integrals on portions of $\Gamma$; next we note that ${\mathcal M} u + f=0$ in every $T \in {\mathcal T}_h \setminus {\mathcal S}_h$; this is also 
true of elements $T$ not belonging to ${\mathcal Q}_h$. 
Finally we recall that $ {\mathcal M} \tilde{u}  + f$ vanishes identically in $T^{'} (= T \cap \Omega$) $\forall T \in {\mathcal Q}_h$. In short we can write:  
\begin{equation}
\label{ahFhprime}
|a_h(\tilde{u},v)-F_h(v)| = \displaystyle \sum_{T \in {\mathcal Q}_h} \left| \int_{\Delta_T} v {\mathcal M} \tilde{u} \right| 
\leq \displaystyle \sum_{T \in {\mathcal Q}_h} \parallel {\mathcal M} \tilde{u} \parallel_{0,\Delta_T} \parallel v \parallel_{0,\Delta_T} \!.
\end{equation}   
Now taking into account that $v \equiv 0$ on $\Gamma_h$ and recalling the constant $C_{\Gamma}$ defined in Proposition \ref{prop01}, it holds :   
$|v({\bf x})| \leq C_{\Gamma} h_T^2 \parallel {\bf grad}\; v \parallel_{0,\infty,\Delta_T}$, $\forall {\bf x} \in \Delta_T$. 
Therefore using \eqref{L2TDelta} we have $\parallel {\bf grad} \; v \parallel_{0,\infty,\Delta_T} \leq {\mathcal C}_J h_T^{-1} 
\parallel {\bf grad} \; v \parallel_{0,T}$. Noticing that the measure of $\Delta_T$ is bounded above by $C_{\Gamma} h_T^3$, after straightforward calculations we obtain for a certain mesh-independent constant $C_R$:
\begin{equation}
\label{fifthbound}
\parallel {\mathcal M} \tilde{u} \parallel_{0,\Delta_T} \parallel v \parallel_{0,\Delta_T} \leq C_R h_T^{5/2} \parallel  {\mathcal M} \tilde{u}  \parallel_{0,\Delta_T} \parallel 
{\bf grad} \; v \parallel_{0,T} \; \forall T \in {\mathcal Q}_h.
\end{equation}
Now combining (\ref{fifthbound}) into (\ref{ahFhprime}) and applying the Cauchy-Schwarz inequality, we easily come up with,
\begin{equation}
\label{sixthbound}
|a_h(\tilde{u},v)-F_h(v)| \leq C_R h^{5/2} \displaystyle \parallel {\mathcal M} \tilde{u} \parallel_{0,\tilde{\Omega}} \parallel {\bf grad} \; v \parallel_{0,h}.
\end{equation}
Finally plugging (\ref{sixthbound}) into (\ref{fourthbound}) and taking into account \eqref{interP2} we easily establish the validity of error estimate (\ref{estimateP2}). \QED \\ 

\begin{theorem}
\label{P3} 
Let $k=3$ and assume that $u \in H^4(\Omega)$. Provided $h$ is sufficiently small, there exists a mesh-independent constant $C_3$ such that the unique solution $u_h$ to (\ref{Poissonh}) satisfies:
\begin{equation}
\label{estimateP3} 
\parallel {\bf grad}(u - u_h) \parallel_{0,h}^{'} \leq C_3 h^{3} [| \tilde{u} |_{4,\tilde{\Omega}} + h^{1/2} \parallel {\mathcal M} \tilde{u}  \parallel_{0,\infty,\tilde{\Omega}}]
\end{equation}
where $\tilde{u} \in H^4(\tilde{\Omega})$ is the regular extension of $u$ to $\tilde{\Omega}$ constructed in accordance to Stein et al. \cite{Stein}.
\end{theorem}

\prov First of all we point out that, according to the Sobolev Embedding Theorem \cite{Adams}, ${\mathcal M} \tilde{u} \in L^{\infty}(\tilde{\Omega})$, since $\tilde{u} \in 
H^4(\tilde{\Omega})$ by assumption. Furthermore taking $w=I_h(\tilde{u})$ and using \eqref{interpolerrorh} we obtain,
\begin{equation}
\label{interP3}
\parallel {\bf grad}(\tilde{u} - w) \parallel_{0,h} \leq \bar{C}_{\Omega} h^3|\tilde{u}|_{4,\tilde{\Omega}}.
\end{equation}  
Now following the same steps as in the proof of Theorem \ref{P2} up to equation (\ref{ahFhprime}), the latter becomes:
\begin{equation}
\label{ahFh3}
|a_h(\tilde{u},v)-F_h(v)| \leq C_{\Gamma} \displaystyle \sum_{T \in {\mathcal Q}_h}  h_T^{3} \parallel {\mathcal M} \tilde{u} \parallel_{0,\infty,\Delta_T} 
\parallel v \parallel_{0,\infty,\Delta_T}, 
\end{equation}
Akin to the previous proof, using \eqref{L2TDelta} we note that, 
\begin{equation}
\label{invers} \parallel v \parallel_{0,\infty,\Delta_T} \leq C_{\Gamma} h_T^2 \parallel {\bf grad}\; v \parallel_{0,\infty,\Delta_T} 
\leq C_{\Gamma} {\mathcal C}_J h_T \parallel {\bf grad}\; v \parallel_{0,T}. 
\end{equation}
Combining (\ref{ahFh3}) with (\ref{invers}) we come up with,
\begin{equation}
\label{ahFh4}
|a_h(\tilde{u},v)-F_h(v)| \leq C_{\Gamma}^2 {\mathcal C}_J \parallel {\mathcal M} \tilde{u} \parallel_{0,\infty,\tilde{\Omega}} \displaystyle \sum_{T \in {\mathcal Q}_h} h_T^4 
\parallel {\bf grad} \; v \parallel_{0,T},  
\end{equation}
Further applying the Cauchy-Schwarz inequality to the right hand side of (\ref{ahFh4}) we easily obtain:
\begin{equation}
\label{ahFh5}
|a_h(\tilde{u},v)-F_h(v)| \leq C_{\Gamma}^2 {\mathcal C}_J h^{7/2} \parallel {\mathcal M} \tilde{u} \parallel_{0,\infty,\tilde{\Omega}} 
\left[ \sum_{T \in {\mathcal Q}_h} h_T \right]^{1/2} 
\parallel {\bf grad} \; v  \parallel_{0,h}.
\end{equation}
From the assumptions on the mesh for a certain mesh-independent constant $C(\Gamma)$ it holds: 
\begin{equation}
\label{CGamma} 
\sum_{T \in {\mathcal Q}_h} h_T  \leq C(\Gamma).
\end{equation}  
Plugging \eqref{CGamma} into (\ref{ahFh5}) and the resulting relation into (\ref{fourthbound}), the result easily follows,  
taking into account \eqref{interP3}. \QED \\

Akin to Theorem \ref{theorem1bis}, it is possible to establish error estimates in the $L^2$-norm in the case of a non convex $\Omega$, by requiring some more regularity from the solution $u$ of \eqref{Poisson}. However,  unless the assumptions of Theorem \ref{theorem2} hold, optimality is not attained for $k>2$. This is because of the absence of $u$ from the non-empty domain $\tilde{\Delta}_h : = \Omega_h \setminus \Omega$, whose area is an invariant $O(h^2)$ whatever $k$. Roughly speaking, integrals in $\tilde{\Delta}_h$ of expressions in terms of the solution $u_h$ of \eqref{Poissonh} dominate the error, in such a way that those terms cannot be reduced to less than an $O(h^{7/2})$, even under additional regularity assumptions. \\
\indent Most steps in the proof of the following result rely on arguments essentially identical to those already exploited to prove Theorem \ref{theorem1bis}. Therefore we will focus on aspects specific to the non convex case.   
 
\begin{theorem}
\label{L2P2} Let $k=2$. Assume that $\Omega$ is not convex and $u \in H^{3+r}(\Omega)$ for 
$r=1/2+\epsilon$, $\epsilon>0 $ being arbitrarily small. Then provided $h$ is sufficiently small the following error estimate holds:
\begin{equation}
\label{L2estP2} \parallel u - u_h \parallel_{0,h}^{'} \leq \tilde{C}_{0} h^{3} [G(\tilde{u})
+ \parallel \tilde{u} \parallel_{3+r,\tilde{\Omega}} ],   
\end{equation} 
where $\tilde{C}_{0}$ is a mesh-independent constant and $G(\tilde{u}):= | \tilde{u} |_{3,\tilde{\Omega}} 
+ h^{1/2} \parallel \nu \Delta \tilde{u} - {\bf b} \cdot {\bf grad} \; \tilde{u} \parallel_{0,\tilde{\Omega}}$.\\
\end{theorem}

\prov Let $\bar{u}_h$ be the function defined in $\Omega$ by $\bar{u}_h:=u_h-u$. $v \in H^1_0(\Omega)$ being the function satisfying \eqref{adjoint}-\eqref{adjoint1}, we have: 
\begin{equation}
\label{L2estP2bis}
\parallel \bar{u}_h \parallel_{0,h}^{'} \leq \parallel \bar{u}_h \parallel_{0} \leq 
C(\Omega) \displaystyle \frac{- \int_{\Omega} \bar{u}_h(\nu \Delta v + {\bf b} \cdot {\bf grad} \; v)}{\parallel v \parallel_{2}}.
\end{equation}
Now we recall the set $\Gamma_h^{'} = \Omega_h \cap \Gamma$ and note that $length(\Gamma_h^{'})>0$. Thus using integration by parts we easily obtain,  
\begin{equation}
\label{L2estP2ter}
\parallel \bar{u}_h \parallel_{0,h}^{'} \leq 
C(\Omega) \displaystyle \frac{b_{1h}(\bar{u}_h,v) + a^{'}_h(\bar{u}_h,v)+a_{\Delta_h}(\bar{u}_h,v)}{\parallel v \parallel_{2}},
\end{equation}
where the bilinear forms $b_{1h}$ and $a_{\Delta_h}$ are defined in \eqref{b1h} and \eqref{aDeltah}, respectively, and  
\begin{equation}
\label{tildeah}
a^{'}_h(w,z) 
:= \int_{\Omega^{'}_h} 
[\nu {\bf grad} \; w \cdot {\bf grad} \; z + ({\bf b} \cdot {\bf grad}\; w)z] \mbox{ for } w,z \in H^1(\Omega).
\end{equation}
On the other hand, since $f = 0$ in $\Omega_h \setminus \Omega$, $\forall v_h \in V_h$ we have,
\begin{equation}
\label{L2estP2qua}
a_h(u_h,v_h)= \int_{\Omega_h^{'}} [-\nu \Delta u + {\bf b} \cdot {\bf grad} \;u]v_h =  
- \displaystyle \nu \int_{\Gamma^{'}_h} \frac{\partial u}{\partial n}v_h + a^{'}_h(u,v_h).
\end{equation}
Recalling the definition of ${\mathcal Q}_h$ in Subsection 2.1 and noting that here 
 $\Delta_T$ is the interior of the set $T \setminus \Omega$ for every $T \in {\mathcal Q}_h$, we define 
\begin{equation}
\label{b5h}
b_{5h}(w,z) := \displaystyle \sum_{T \in {\mathcal Q}_h} \int_{\Delta_T} [-\nu \Delta w + {\bf b} \cdot {\bf grad}\; w]z, 
\forall w \in W_h \mbox{ and } \forall z \in V_h.
\end{equation}
Denoting also by $\partial T$ the set $\Gamma \cap T$ we further set,
\begin{equation}
\label{b6h}
b_{6h}(w,z) := \nu \displaystyle \sum_{T \in {\bf Q}_h} 
\int_{\partial T} \frac{\partial w}{\partial n}z \; \forall  
w \in W_h \cup H^2(\Omega) \mbox{ and } z \in V_h.
\end{equation}
It easily follows from \eqref{L2estP2qua} that 
\begin{equation}
\label{L2estP2qui}
-a^{'}_h(\bar{u}_h,v_h) + b_{5h}(u_h,v_h) + b_{6h}(\bar{u}_h,v_h)=0 \; \forall v_h \in V_h.
\end{equation}
Taking $v_h = \Pi_h(v)$, recalling that $e_h(v):=v-\Pi_h(v)$ and plugging \eqref{L2estP2qui} into \eqref{L2estP2ter} we come up with,
\begin{equation}
\label{L2estP2sex}
\parallel \bar{u}_h \parallel_{0,\tilde{\Omega}_h} \leq 
C(\Omega) \displaystyle \frac{b_{1h}(\bar{u}_h,v) + 
b_{5h}(u_h,\Pi_h(v))+b_{6h}(\bar{u}_h,\Pi_h(v))+a^{'}_h(\bar{u}_h,e_h(v))+a_{\Delta_h}(\bar{u}_h,v)}{\parallel v \parallel_{2}}.
\end{equation}
On the other hand, recalling $b_{2h}$ given by \eqref{b2h} and using integration by parts we have
\begin{equation}
\label{L2estP2sept}
a_{\Delta_h}(\bar{u}_h,v)=a_{\Delta_h}(\bar{u}_h,e_h(v)) +  
 \displaystyle \sum_{T \in {\mathcal S}_h \setminus {\mathcal Q}_h} \int_{\partial T} \displaystyle 
\nu \frac{\partial \bar{u}_h}{\partial n} \Pi_h(v) + b_{2h}(\bar{u}_h,\Pi_h(v)).
\end{equation} 
Thus recalling $b_{3h}$ and $b_{4h}$ respectively defined by \eqref{b3h} and \eqref{b4h}, we finally obtain:
\begin{equation}
\label{L2estP2oct}
\left\{
\begin{array}{l}
\parallel \bar{u}_h \parallel_{0,h}^{'} \leq 
C(\Omega) \displaystyle \frac{L(\bar{u}_h,v)+b_{5h}(u_h,v_h)+a^{'}_h(\bar{u}_h,e_h(v))}{\parallel v \parallel_{2}},\\
\mbox{where }\\
L(\bar{u}_h,v):=b_{1h}(\bar{u}_h,v) + b_{2h}(\bar{u}_h,\Pi_h(v)) + b_{3h}(\bar{u}_h,\Pi_h(v))+ 
b_{4h}(\bar{u}_h,e_h(v)).
\end{array}
\right.
\end{equation}
The estimation of $a^{'}_h(\bar{u}_h,e_h(v))$ is a trivial variant of 
\eqref{L2est3bis}, that is,
\begin{equation}
\label{estildeah}
a^{'}_h(\bar{u}_h,e_h(v)) \leq C_2 C^{'}_V \tilde{A} h^{3} G(\tilde{u}) |v|_{2},
\end{equation}
where $C^{'}_V$ is an interpolation error constant such that 
\begin{equation}
\label{intildeP1}
\parallel {\bf grad} [v-\Pi_h(v)]\parallel_{0,h}^{'} \leq C^{'}_V h | v |_{2}.
\end{equation}
\indent The bilinear forms $b_{ih}$, $i=1,2,3,4$ were studied in Theorem \ref{theorem1bis}. The corresponding estimates here are qualitatively the same taking $k=2$, if we replace here and there $| u |_{3}$ 
by $G(\tilde{u})$. Hence all that is left to do is to estimate $b_{5h}(u_h,v_h)$. With this aim we proceed as follows:\\

Since $|v_h({\bf x})| \leq C_{\Gamma} h_T^2 \parallel {\bf grad}\; v_h \parallel_{0,\infty,T}$ $\forall {\bf x} \in \Delta_T$ and $\forall T \in {\mathcal Q}_h$ for $v_h \in V_h$, by a straightforward argument we can write
\begin{equation}
\label{step1}
b_{5h}(u_h,\Pi_h(v)) \leq \displaystyle \sum_{T \in {\mathcal Q}_h} C^2_{\Gamma} h_T^5 
[ \nu  \parallel \Delta u_h \parallel_{0,\infty,T} +
\parallel {\bf b} \parallel_{0,\infty} \parallel {\bf grad} \; u_h \parallel_{0,\infty,T} ]
\parallel {\bf grad}\; \Pi_h(v) \parallel_{0,\infty,T}.
\end{equation}
Since all the components of $[{\bf grad}\;\Pi_h(v)]_{|T}$ and $[H(u_h)]_{|T}$ are in ${\mathcal P}_0(T)$ and those of $[{\bf grad}\;u_{h}]_{|T}$ are in ${\mathcal P}_1(T)$, in all the norms involving $\Pi_h(v)$ and $u_h$ in  
\eqref{step1}, $T$ can be replaced by $T^{'}$. Moreover we have,
\begin{equation}
\label{step2}
\begin{array}{l}
\nu  \parallel \Delta u_h \parallel_{0,\infty,T^{'}} + \parallel {\bf b} \parallel_{0,\infty} \parallel {\bf grad} \; u_h \parallel_{0,\infty,T^{'}} 
\leq \nu  \parallel \Delta u \parallel_{0,\infty,T^{'}} + \parallel {\bf b} \parallel_{0,\infty} \parallel {\bf grad} \; u \parallel_{0,\infty,T^{'}} \\
+\nu  \parallel \Delta \bar{u}_h \parallel_{0,\infty,T^{'}} + \parallel {\bf b} \parallel_{0,\infty} \parallel {\bf grad} \; \bar{u}_h \parallel_{0,\infty,T^{'}} 
\end{array}
\end{equation}
After plugging \eqref{step2} into \eqref{step1} we use \eqref{L2TDelta} with $w=\Pi_h(v)$ and $k=1$ together with \eqref{Djwh} with $w=\tilde{u}$ and $w_h=u_h$ to obtain , 
\begin{equation}
\label{step3}
\begin{array}{l}
b_{5h}(u_h,\Pi_h(v)) \leq C^2_{\Gamma} {\mathcal C}_J \displaystyle \sum_{T \in {\mathcal Q}_h}  h_T^4 \parallel {\bf grad}\; \Pi_h(v) \parallel_{0,T^{'}} [
(\nu \parallel \Delta u \parallel_{0,\infty} + \parallel {\bf b} \parallel_{0,\infty}\parallel {\bf grad} \; u \parallel_{0,\infty} )
\\
 +( \nu \sqrt{2} {\mathcal C}_2 h_T^{-2} + \parallel {\bf b} \parallel_{0,\infty} {\mathcal C}_1 h_T^{-1}) ( 
\parallel {\bf grad} \; \bar{u}_h \parallel_{0,T^{'}} + h_T^2 | u |_{3,T{'}} + h_T^{2+r} \parallel \tilde{u} \parallel_{3+r,T}) ].
\end{array}
\end{equation} 
 Therefore, applying the Cauchy-Schwarz inequality and recalling \eqref{CGamma}, we conclude that there a mesh-independent constant $\tilde{C}_5$ such that, 
\begin{equation}
\label{step4}
\begin{array}{l}
b_{5h}(u_h,\Pi_h(v)) \leq \tilde{C}_5 \parallel {\bf grad}\; \Pi_h(v) \parallel_{0,h}^{'} [ h^{7/2} \parallel u \parallel_{2,\infty} \\ 
+ h^2 (\parallel {\bf grad} \; \bar{u}_h \parallel_{0,h}^{'} +  h^2| u |_{3} + h^{2+r} \parallel \tilde{u} \parallel_{3+r,\tilde{\Omega}}) ].
\end{array}
\end{equation}
Now since $h<1$, recalling \eqref{estimateP2} and noting that by the Sobolev embedding Theorem (cf. \cite{Adams}) there exists a constant $C_e$ depending only on $\Omega$ such that $\parallel u \parallel_{2,\infty} \leq C_e  \parallel u \parallel_{3+r}$, 
from \eqref{step4} we obtain:
\begin{equation}
\label{step5}
b_{5h}(u_h,\Pi_h(v)) \leq \bar{C}_5  [ h^{7/2} \parallel \tilde{u} \parallel_{3+r,\tilde{\Omega}} 
 + h^4 G(\tilde{u}) ] \parallel {\bf grad}\; \Pi_h(v) \parallel_{0,h}^{'}.
\end{equation}
where $\bar{C}_5$ is a suitable mesh-independent constant. \\
On the other hand from \eqref{intildeP1} we easily derive, 
\begin{equation}
\label{interpol}
\parallel {\bf grad}\; \Pi_h(v) \parallel^{'}_{0,h} \leq C_{\Pi} \parallel v \parallel_{2}.
\end{equation} 
with $C_{\Pi}=\sqrt{1+(C^{'}_V)^2 h_0^2}$. Plugging \eqref{interpol} into \eqref{step5}, setting 
$C_{b5}:= C_{\Pi} \bar{C}_5$, we conclude that,
\begin{equation}
\label{estimateb5}
b_{5h}(u_h,v_h) \leq  C_{b5} h^{7/2} \displaystyle \left\{ h^{1/2} G(\tilde{u}) 
+ \parallel \tilde{u} \parallel_{3+r,\tilde{\Omega}} \right\} \parallel v  \parallel_{2}.
\end{equation} 

\noindent Finally recalling \eqref{L2estP2oct} together with \eqref{estildeah}, \eqref{estimateb1}, \eqref{estimateb2}, \eqref{estimateb3} and \eqref{estimateb4}, estimate \eqref{estimateb5} completes the proof. \QED \\ 
 
\begin{remark} 
Several computations with the method studied in this section were carried out. 
It comes out from the thus generated numerical results reported in about ten publications so far, that its performance is completely satisfactory. This method was even found to be superior to the isoparametric technique in terms of accuracy, while both methods showed to be fairly equivalent as far as CPU time is concerned. 
For further details on this numerical experimentation the author refers to 
\cite{CFM2017}, \cite{arXiv2D}, \cite{Maugin} and \cite{CAMWA}. Numerical experiments in three-dimensional space can be found in \cite{arXiv3D} and \cite{IMAJNA}. A sample of those results is supplied in the Appendix hereafter. \QED
\end{remark}

\section{Hermite elements with normal-derivative degrees of freedom}

The technique described in Section 2 for a simple model problem is an effective tool to handle finite-element degrees of freedom of various types to be prescribed on curvilinear boundaries. This is because 
it extends in a rather straightforward manner to classes of finite-element methods other than the classical Lagrange family, as seen in \cite{AMIS} and \cite{CFM2017}. 
In this section we illustrate this assertion by considering the case of the biharmonic problem in a smooth two-dimensional domain with Dirichlet boundary conditions. The equation under consideration governs the deflexion $\psi$ of a thin plate with a clamped curved edge 
$\Gamma$, that is the boundary of the domain $\Omega$ of the plane schematically occupied by the plate, under the action of a distribution of forces $\tilde{f}$ perpendicular to this plane. In a suitable system of units the problem to solve is 
\begin{equation}
\label{psi}
\left\{
\begin{array}{ll}
\Delta^2 \psi = f & \mbox{ in } \Omega\\
\psi = 0 \mbox{ and } \displaystyle \frac{\partial \psi}{\partial n} = 0 & \mbox{ on } \Gamma,
\end{array}
\right.
\end{equation}
where $f$ is proportional to $\tilde{f}$ and supposedly belongs to $L^2(\Omega)$.\\ 
Equation \eqref{psi} can be rewritten in the following equivalent variational form:
\begin{equation}
\label{Weakpsi}
\left\{
\begin{array}{l}
\mbox{Find } \psi \in H^2_0(\Omega) \mbox{ such that }  
c(\psi,\phi) = F(\phi) \; \forall \phi \in H^2_0(\Omega),\\
\mbox{where,}  \\
c(\psi,\phi) :=  \int_{\Omega} \Delta \psi \Delta \phi \mbox{ and } F(\phi) := \int_{\Omega} f \phi. 
\end{array}
\right.
\end{equation}
The bilinear form $c$ is coercive over $[H^2_0(\Omega)]^2$ equipped with the norm $| \cdot |_{2}$, i.e., the usual semi-norm of $H^2(\Omega)$. Indeed  by density arguments and integration by parts we easily conclude that $c(\phi,\phi) = \int_{\Omega} | \Delta \phi |^2 =  \int_{\Omega} | H (\phi )|^2 = | \phi |_{2}^2$.\\
Hereafter we use the notations related to a triangular mesh introduced in Section 2. \\
As it is well-known, in principle, conforming finite-element methods to solve a fourth-order boundary-value problem must be based on test- and trial functions of the $C^1$-class. Most popular methods falling into this  category for triangular meshes, have first-order normal-derivative degrees of freedom at points in the interior of the edges such as mid-points. In order to concentrate on the main aspects of our technique's application to this kind of degree of freedom, we restrict the presentation to the case of the Clough-Tocher element (cf. 
\cite{C-T}, also known as the Hsieh-Clough-Tocher element \cite{Ciarlet}). With minor modifications the analysis extend to other $C^1$ finite element methods with this type of degree of freedom, such as the conforming Zienkiewicz elements among other methods listed in \cite{Ciarlet}. We recall that the degrees of freedom for the Clough-Tocher element are the values of the function and its first order derivatives at the vertexes, and the first order normal derivative at the edge mid-points. This set of twelve degrees of freedom for each triangle is sufficient to uniquely define a piecewise cubic function of the $C^1$-class therein, whose precise construction is described in \cite{Ciarlet}. Coincidence of these degrees of freedom at inter-element boundaries also ensure the continuity along the edges of the mesh of both the resulting piecewise cubic function and its first order derivative in the direction normal to the edges. Moreover, as long as all the degrees of freedom attached to the nodes located on $\Gamma_h$ vanish, the corresponding function and its (normal) derivative in the direction of ${\bf n}_h$ will vanish on $\Gamma_h$. We denote by $\Phi_h$ the space of functions satisfying such conditions. \\
We also recall that, in the case of a polygonal domain, and as long as the solution belongs to $H^4(\Omega)$, the Clough-Tocher method is of the second order in the standard (semi-)norm of $H^2(\Omega)$. \\
Extending $f$ by zero in $\Omega_h \setminus \Omega$ and still denoting the resulting function by $f$, let us first 
approximate \eqref{Weakpsi} by 
\begin{equation}
\label{Weakpsih}
\left\{
\begin{array}{l}
\mbox{Find } \psi^{'}_h \in \Phi_h \mbox{ such that } c_h(\psi^{'}_h,\phi) = F_h(\phi) \; \forall \phi \in \Phi_h,\\
\mbox{where,}   \\
c_h(\psi,\phi) :=  \int_{\Omega_h} \Delta \psi \Delta \phi \mbox{ and } F_h(\phi) := \int_{\Omega_h} f \phi. 
\end{array}
\right.
\end{equation} 
It is not difficult to verify that an order reduction from $2$ to $3/2$ in the above sense will occur, if we prescribe 
a zero value to the solution's normal derivative at the mid-points of the edges contained in $\Gamma_h$. \\
In order to overcome this issue it is necessary to prescribe normal derivative degrees of freedom at locations on the true boundary. As far as the author can see, a true isoparametric version of the Clough-Tocher triangle cannot be worked out. In \cite{Mansfield} a quadratic sub-parametric version for curved domains was studied. However this technique has the drawback of transforming method's original degrees of freedom into non trivial combinations thereof. That is why we propose instead a technique entirely analogous to the modification of $V_h$ into $W_h$ described in Section 2. \\
\indent In order to avoid fastidious calculations, in the remainder of this section we voluntarily skip some details of  our analysis based on arguments analogous to those already exploited in Sections 2 and 3. \\ 
\indent Here we replace $\Phi_h$ as a trial-space by a space $\Psi_h$ defined in the same manner as the former, except for the fact that the degrees of freedom corresponding to a normal derivative at the mid-point $M$ of each edge contained in $\Gamma_h$ is replaced with the (zero) value of the normal derivative at a point $P \in \Gamma$ close to $M$. Such a point $P$ can be for instance the closest intersection with $\Gamma$ of the perpendicular to this edge passing through $M$, as indicated in Figure 3. In the same figure a self-explanatory illustration of the degrees of freedom for the Clough-Tocher element related to a triangle in $S \in {\mathcal S}_h$ is 
provided. ${\bf n}$ denotes the unit outer normal vector along $\Gamma$ at point $P$ and the corresponding 
arrow also represents the normal-derivative degree of freedom attached to $P$. \\

\begin{figure}[H]
\label{fig3}
\centerline{\includegraphics[width=3.5in]{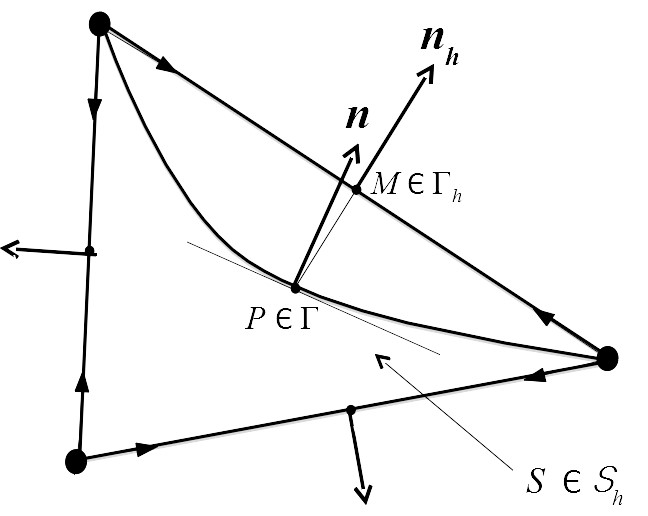}}
\vspace*{8pt}
\caption{Degrees of freedom for space $\Psi_h$ restricted to a triangle $S$ next to a concave portion of $\Gamma$}
\end{figure}

The construction of $\Psi_h$ as above is feasible according to,
\begin{lemma}
\label{lemmaCT}
Referring to \cite{Ciarlet} for further details, let $S$ be a triangle in ${\mathcal S}_h$ and $S_1,S_2,S_3$ be the three disjoint triangles whose union is $S$, having the centroid of $S$ as a vertex.   
Let $P_{CT}(S)$ be the space of functions of the $C^1$-class in $S$, whose restriction to each $S_i$, $i=1,2,3$ is a cubic function, defined by the following set of degrees of freedom denoted by ${\mathcal F}_i, \; i=1,2,\ldots,12$: the values at the vertexes of $S$, the values of the first order derivatives at the vertexes of $S$, the values of the outer normal derivatives at the mid-point of the edges not contained in $\Gamma_h$, 
and the value of $\partial (\cdot)/\partial n$ at the point $P \in \Gamma$ as indicated in Figure 3.      
Then given a set of $12$ real numbers $b_{i}$, $i=1,2,\ldots,12$ there exists a unique function $w \in P_{CT}(S)$ such that ${\mathcal F}_i(w)=b_i$ for $i=1,2,\ldots,12$. 
\end{lemma}   
 
\prov Let $S_1$ be the sub-triangle of $S$ having an edge $e_S$ contained in $\Gamma_h$. For convenience we let the normal derivative at the point $P \in \Gamma$ associated with the mid-point $M$ of $e_S$ to be ${\mathcal F}_1$.\\
Let also $\{{\mathcal G}_j\}_{j=1}^{12}$ be the set of degrees of freedom of the standard Clough-Tocher element 
and $\{\zeta_j\}^{12}_{j=1}$ be corresponding set of canonical basis functions. Clearly enough, ${\mathcal G}_j=
{\mathcal F}_j$ for $j=2,\ldots,12$ and ${\mathcal G}_1(\cdot):=[\partial (\cdot)/\partial n_h](M)$. Similarly to Lemma \ref{lemma1} the key to the problem is the estimation of the difference between the $12 \times 12$ identity matrix ${\mathcal I}$ and the coefficients of the $12 \times 12$ matrix $K=\{k_{ij}\}$ whose entry $k_{ij}$ is given by ${\mathcal F}_{i}(\zeta_j)$. Letting $E =\{e_{ij}\}$ be the matrix ${\mathcal I}-K$ we easily find out that $e_{ij}=0$ for every $i>1$ and $e_{1j}=[{\mathcal F}_1-{\mathcal G}_1]({\zeta_j})$ for $j=1,2,\ldots,12$.  It follows that the matrix $K$ is a lower triangular matrix, and thus its determinant equals  
$1+ e_{11}$ with $e_{11}= [{\mathcal F}_1-{\mathcal G}_1]({\zeta_1})$. \\
Now we observe that $\zeta_1$ vanishes along $e_S$, which implies in particular that its tangential derivative vanishes at $M$. This means that the partial derivative of $\zeta_1$ in the direction of ${\bf n}$ at $M$ equals the partial derivative of $\zeta_1$ in the direction of ${\bf n}_h$ at the same point multiplied by ${\bf n} \cdot {\bf n}_h$. Denoting by $\theta$ be angle between ${\bf n}$ and ${\bf n}_h$ we thus have:
\begin{equation}
\label{e111} 
e_{11} = \displaystyle \left[\frac{\partial \zeta_1}{\partial n} \right](P)-\left[\frac{\partial \zeta_1}{\partial n} \right](M) - 
[(cos \theta)^{-1} - 1] \left[\frac{\partial \zeta_1}{\partial n} \right](M). 
\end{equation}
Recalling the function $f_S$ expressing the boundary $\Gamma$ at the level of $S$ we note that $(cos \theta)^{-1} = \sqrt{1+[tan \theta]^2}$ 
$= \sqrt{1+[f_S^{'}(P)]^2}$. Thus after straightforward calculations $(cos \theta)^{-1}-1$ is found to be bounded above by $C_{B} h_S$, 
with $C_B = {\mathcal H}_{max}/\sqrt{2}$.  \\
Now referring to the notations employed in Lemma \ref{lemma1}, we readily have:
\begin{equation}
\label{e112} 
|e_{11}| \leq \displaystyle \left \{ C_{\Gamma}  \max_{Q \in \overline{PM}} \left| \left[ \frac{\partial^2 \zeta_1}{\partial n \partial n_h}\right](Q)\right|h_S^2 +  C_{B} \displaystyle \left| \left[ \frac{\partial \zeta_1}{\partial n}\right](M) \right| h_S \right\} 
\end{equation}  
Since the first order derivatives of $\zeta_1$ in $S$ are bounded by a constant independent of $S$ and the second order derivatives of $\zeta_1$ in $S$ can be bounded by $h_S^{-1}$ multiplied by another constant 
independent of $S$, $det K= 1+e_{11} \neq 0$ if $h_S$ is sufficiently small. \QED \\
     
Now quite naturally, we modify the approximate problem into, 
\begin{equation}
\label{Weakpsihn}
\left\{
\begin{array}{l}
\mbox{Find } \psi_h \in \Psi_h \mbox{ such that,}  \\
c_h(\psi_h,\phi) = F_h(\phi) \; \forall \phi \in \Phi_h,
\end{array}
\right.
\end{equation} 
We equip both $\Psi_h$ and $\Phi_h$ with the standard semi-norm of $H^2(\Omega_h)$ 
denoted in the sequel by $| \cdot |_{2,h}$, or otherwise stated, with the norm $\parallel H(\cdot) \parallel_{0,h}$. Thanks to the fact that every function $\eta$ belonging to either $\Psi_h$ or $\Phi_h$  
together with its gradient, vanishes at the vertexes of $\Gamma_h$, $| \eta |_{2,h}$ is a norm indeed. \\
We next prove the counterpart of Proposition \ref{propo1} for problem \eqref{Weakpsihn}. We omit details whenever the underlying argument can be found in the proof thereof.
\begin{e-proposition}
\label{propo2}
Provided $h$ is sufficiently small problem \eqref{Weakpsihn} has a unique solution. Moreover there exists a constant $\gamma > 0$ independent of $h$ such that,
\begin{equation}
\label{stab}
\forall \psi \in \Psi_h \neq 0, \displaystyle \sup_{\phi \in \Phi_h \setminus \{ 0 \}} \frac{c_h(\psi,\phi)}{| \psi |_{2,h} | \phi |_{2,h}} 
\geq \gamma.
\end{equation}  
\end{e-proposition}

\prov First of all we observe that 
\begin{equation}
\label{equivalence}
 \int_{\Omega_h} \Delta \psi \Delta \phi = \int_{\Omega_h} H(\psi):H(\phi) \; \forall \psi \in \Psi_h \mbox{ and } \forall \phi \in \Phi_h. 
\end{equation}
\eqref{equivalence} is a straightforward consequence of the application of First Green's formula in all the sub-triangles $S_1,S_2,S_3$ of the mesh back and forth. To check this identity it suffices to take into account the commutativity of the operators ${\bf grad}$ and $\Delta$, the continuity of both $\psi$ and $\phi$ and also of their normal derivatives along the inner edges, and the fact that 
both $\phi$ and its normal derivative vanish on $\Gamma_h$. In view of this the summation of all the jump terms along the interfaces of elements and sub-elements cancel out, and all the integrals along edges of $\Gamma_h$ vanish.\\
Now given $\psi \in \Psi_h$, we construct an associated $\phi \in \Phi_h$ whose degrees of freedom coincide with those of $\psi$, except for the outer normal derivative at the mid-points $M$ of edges $e_S \subset \Gamma_h$ for every $S \in {\mathcal S}_h$. By construction the normal derivatives of $\phi$ vanish at all those points.\\ 
Then, taking into account \eqref{equivalence} we have,  
\begin{equation}
\label{chpsiphi}
\begin{array}{l} 
c_h(\psi,\phi) = \displaystyle \sum_{T \in {\mathcal T}_h} \int_T  |H(\psi)|^2 - \displaystyle \sum_{S \in {\mathcal S}_h} \int_S 
 H(\psi) : H(\chi_S(w)),
\end{array} 
\end{equation}
\noindent where $\chi_S(\psi) = \displaystyle \left[ \frac{\partial \psi}{\partial n_h} \right](M) \zeta_M$, $\zeta_M$ being the basis function $\zeta_1$ specified in Lemma \ref{lemmaCT} for triangle $S$. \\
\noindent It is noticeable that $\max_{Q \in S}|\zeta_M(Q)|$ is an $O(h_S)$. Thus from standard results it holds for a mesh-independent constant $C_{\zeta}$:
\begin{equation}
\label{zetaestim}
\parallel H(\zeta_M) \parallel_{0,S} \leq C_{\zeta}.
\end{equation} 
On the other hand, since $\partial \psi/\partial n$ equals zero at $P \in \Gamma$ by construction, using the very same arguments as those leading to \eqref{e111} and \eqref{e112} we have:
\begin{equation}
\label{normderiv} 
\displaystyle \left| \left[ \frac{\partial \psi}{\partial n_h} \right](M) \right| = 
(cos \theta)^{-1} \displaystyle \left| \left[ \frac{\partial \psi}{\partial n} \right](M) \right|
\leq (1+C_B h_S) C_{\Gamma} h_S^2 \displaystyle \max_{Q \in \overline{PM}} \left| \left[ \frac{\partial^2 \psi}{\partial n \partial n_h}\right](Q) \right|. 
\end{equation}
Moreover, from an inverse inequality analogous to \eqref{invinity}, we may write for a suitable constant $C_{\infty}$:
\begin{equation}
\label{inversineq}
\parallel  H(\psi) \parallel_{0,\infty,S} \leq C_{\infty} h_S^{-1} \parallel H(\psi) \parallel_{0,S}.
\end{equation}
Plugging \eqref{zetaestim}, \eqref{normderiv}, and \eqref{inversineq} into (\ref{chpsiphi}), it 
follows that:
\begin{equation}
\label{chpsiphibound} 
c_h(\psi,\phi) \geq \int_{\Omega_h} | H(\psi) |^2 - h \displaystyle \sum_{S \in {\mathcal S}_h} \parallel H(\psi) \parallel_{0,S}^2 
C_{\infty} (1+C_B h) C_{\Gamma} C_{\zeta}   . 
\end{equation}         
From \eqref{chpsiphibound} we readily obtain for a suitable mesh-independent constant $C_H$:
\begin{equation}
\label{chpsiphibelow} 
c_h(\psi,\phi) \geq  (1 - C_H h) |\psi |_{2,h}^2
\end{equation}
On the other hand, similarly to \eqref{normbound} we have  
\begin{equation}
\label{Hessianbound}
 | \phi |_{2,h} \leq | \psi |_{2,h} + | \psi - \phi |_{2,h} \leq (1+ C_H h) 
| \psi |_{2,h}.
\end{equation}
Combining (\ref{chpsiphibelow}) and (\ref{Hessianbound}), provided $h \leq (2 C_H)^{-1}$, we establish (\ref{stab}) with $\gamma = 1/ 3$. \\

\noindent Finally, since $dim(\Phi_h) = dim(\Psi_h)$, \eqref{Weakpsihn} is uniquely solvable. \QED \\

As a consequence of Proposition \ref{propo2}, we can prove an error estimate for problem \eqref{Weakpsihn}. Before doing it we need an approximation result.\\
First we note that, since $H^3(\tilde{\Omega}_h)$ is continuously embedded into $C^1(\overline{\tilde{\Omega}_h})$, $\forall 
\chi \in H^3(\tilde{\Omega}_h)$ we can construct an interpolating function $K_h(\chi) \in \Psi_h$ extended to $\bar{\Omega} \setminus \bar{\Omega}_h$ as 
advocated in the definition of $\Psi_h$, using the degrees of freedom of this space. This is performed by analogy to what we did to define the interpolating function $I_h(w) \in W_h$ of $w \in H^2(\tilde{\Omega}_h)$. Then we have,

\begin{lemma}
For all $\chi \in H^4(\tilde{\Omega})$ there exists a mesh-independent constant $C_X$ such that the interpolating function $K_h(\chi)$ fulfills, 
\begin{equation}
\label{interpCT}
| \chi - K_h(\chi) |_{2,h} \leq C_{X} h^2 | \chi |_{4,\tilde{\Omega}}.
\end{equation}
\end{lemma}

\prov  
Using again the arguments in \cite{BrennerScott}, Subsection 4.4.1, similarly to \eqref{interpolerrorh}, it is readily seen that 
$| \chi - K_h(\chi) |_{2,h} \leq C_{X} h^2 \parallel D^4\chi \parallel_{0,h}$, which implies \eqref{interpCT}. \QED \\

\begin{theorem}
\label{CT}
Provided $h$ is sufficiently small, if the solution $\psi$ of \eqref{psi} is in $H^4(\Omega)$, there exists a mesh-independent constant $C_{CT}$ such that the unique solution $\psi_h$ to (\ref{Weakpsihn}) satisfies:
\begin{equation}
\label{estimCT} 
| H(\psi) - H(\psi_h) |^{'}_{0,h} \leq C_{CT} h^2 J(\tilde{\psi}).
\end{equation}
where $\tilde{\psi}$ is an extension  
of $\psi$ to $\tilde{\Omega}$ belonging to 
$H^4(\tilde{\Omega})$ constructed as advocated in Stein et al. \cite{Stein} and $J(\tilde{\psi}):= | \tilde{\psi} |_{4,\tilde{\Omega}} + h^{5/2} \parallel \Delta^2 \tilde{\psi} \parallel_{0,\tilde{\Omega}}$.
\end{theorem}
\prov 
In this proof we follow the main steps as in Theorems \ref{P2} and \ref{P3}.\\ 
First of all straightforward manipulations lead to, 
\begin{equation}
\label{boundCT1}
| \psi_h - \chi |_{2,h} \leq \displaystyle \frac{1}{\gamma} 
\left[| \tilde{\psi} - \chi |_{2,h} + \displaystyle \sup_{\phi \in \Phi_h \setminus \{0\}} 
\frac{|c_h(\tilde{\psi},\phi)-F_h(\phi)|}{| \phi |_{2,h}} \right] \; \forall \chi \in \Psi_h.
\end{equation}
The term in the numerator of (\ref{boundCT1}) is estimated as follows: Since $\tilde{\psi} \in H^4(\tilde{\Omega})$ and $\Phi_h \subset H^2_0(\Omega_h)$, we apply twice the First Green's formula to $c_h(\tilde{\psi},\phi)$ for $\phi \in \Phi_h$ to obtain,
\begin{equation}
\label{boundCT2}
c_h(\tilde{\psi},\phi) =  \int_{\Omega_h} \Delta^2 \tilde{\psi} \phi.
\end{equation}
Since $\Delta^2 \tilde{\psi} = f$ in $T$ $\forall T \in {\mathcal T}_h \setminus {\mathcal Q}_h$, and also in $T^{'}= T \cap \Omega$ 
$\forall T \in {\mathcal Q}_h$, from \eqref{boundCT2} it follows that,
\begin{equation}
\label{boundCT3}
|c_h(\tilde{\psi},\phi) - F_h(\phi)| = \displaystyle \sum_{T \in {\mathcal Q}_h} \int_{\Delta_T} \Delta^2 \tilde{\psi} \phi.
\end{equation}
or yet,
\begin{equation}
\label{boundCT4}
|c_h(\tilde{\psi},\phi) - F_h(\phi)| \leq \displaystyle \sum_{T \in {\mathcal Q}_h} \parallel \Delta^2 \tilde{\psi} \parallel_{0,\Delta_T} 
\parallel \phi \parallel_{0,\Delta_T}.
\end{equation} 
Now taking into account that both $\phi$ and $\partial \phi/\partial n_h$ vanish identically on $\Gamma_h$ and recalling the constant $C_{\Gamma}$ defined in Proposition \ref{prop01}, by a simple Taylor expansion about points $M \in e_T$ we obtain :   
$|\phi(Q)| \leq C_{\Gamma}^2 h_T^4 \parallel H(\phi)\parallel_{0,\infty,\Delta_T}/2$, $\forall Q \in \Delta_T$. 
Since the area of $\Delta_T$ is bounded above by $C_{\Gamma} h_T^3$, this implies that,
\begin{equation}
\label{boundCT5}
|c_h(\tilde{\psi},\phi) - F_h(\phi)| \leq \displaystyle \sum_{T \in {\mathcal Q}_h} [C_{\Gamma}]^{5/2} h_T^{11/2} \parallel \Delta^2 \tilde{\psi} \parallel_{0,\Delta_T} 
\parallel H(\phi) \parallel_{0,\infty, \Delta_T}/2.
\end{equation} 
Now recalling \eqref{L2TDelta} we have $\parallel H(\phi) \parallel_{0,\infty,\Delta_T} \leq {\mathcal C}_J h_T^{-1} 
\parallel H(\phi) \parallel_{0,T}$. Hence after straightforward calculations we come up with,  
\begin{equation}
\label{boundCT6}
|c_h(\tilde{\psi},\phi) - F_h(\phi)| \leq \displaystyle \sum_{T \in {\mathcal Q}_h} C_{\Psi} h_T^{9/2} \parallel \Delta^2 \tilde{\psi} \parallel_{0,\Delta_T} 
\parallel H(\phi) \parallel_{0,T}, 
\end{equation}
for a certain mesh-independent constant $C_{\Psi}$, or yet, applying the Cauchy-Schwarz inequality, 
\begin{equation}
\label{boundCT7}
|c_h(\tilde{\psi},\phi) - F_h(\phi)| \leq C_{\Psi} h_T^{9/2} \parallel \Delta^2 \tilde{\psi} \parallel_{0,\tilde{\Omega}} 
| \phi |_{2,h}.
\end{equation}

Finally taking $\chi = K_h(\tilde{\psi})$ in \eqref{boundCT1} and plugging (\ref{boundCT7}) and \eqref{interpCT} into the resulting inequality, we easily establish the validity of error estimate (\ref{estimCT}). \QED \\

\section{Possible extensions and conclusions}

To conclude we make some comments on the methodology studied in this work.\\ 


\indent (i) First of all it is important to stress that the assumption on the magnitude of the mesh parameter $h$ made throughout the paper is just a sufficient condition for the formal results given in this work to hold. It is by no means a necessary condition, and actually we can even assert that it is a rather academic hypothesis. 
Indeed good numerical results can be obtained with meshes as coarse as can be, as shown by computations 
such as those reported in \cite{arXiv2D} and \cite{PAMM}. More precisely we have computed repeatedly with meshes for which $h$ was as large as a half diameter of $\Omega$,  
and no problem at all has ever been detected.\\  
\indent (ii) As seen in Section 4, the technique advocated in this work provides a simple and reliable manner to overcome technical difficulties brought about by more complicated problems and interpolations. Actually this feature was  illustrated for the first time in \cite{AMIS}, where our technique was applied to a Hermite analog of the Raviart-Thomas mixed finite element method of the lowest order \cite{RaviartThomas} to solve Maxwell's equations with Dirichlet conditions on the normal derivative. \\
\indent (iii) The solution of boundary value problems with inhomogeneous boundary conditions using our method is straightforward. For instance, in the case of \eqref{Poisson} it suffices to assign the value of $g$ 
at each node belonging to the true boundary $\Gamma$ for any boundary element, that is, any element having an edge contained in $\Gamma_h$. The error estimates derived in this paper trivially extends to this case as the reader can easily figure out. On the other hand in the case of Neumann boundary conditions $\partial u /\partial n = 0$  on $\Gamma$ our method coincides with the standard Lagrange finite element method. Incidentally we recall that in case inhomogeneous Neumann boundary conditions are prescribed, optimality can only be recovered if the linear form $F_h$ is modified in such a way that boundary integrals for elements $T \in {\mathcal S}_h$ are shifted to the curved    
boundary portion of an element $\tilde{T}$ sufficiently close to the one of the corresponding curved element $T^{\Delta}$ or $T^{'}$. But this is an issue that has nothing to do 
with our method, which is basically aimed at resolving those related to the prescription of degrees of freedom for essential boundary conditions. \\
\indent (iv) As the reader has certainly noticed, in order to compute the element matrix and right hand side vector for a triangle in ${\mathcal S}_h$, we have to determine the inverse of an $n_k \times n_k$ matrix. However this extra effort should by no means be a problem at the current state-of-the art of Scientific Computing.\\
\indent (v) It is important to point out that in the case of constant coefficients our method can do without numerical integration to compute element matrices, at least for quadratic and cubic finite elements, as pointed out in Sections 1 and 3. This is another significant advantage thereof over the isoparametric version of the finite element method. Indeed the latter helplessly requires numerical integration for this purpose, since it deals with rational shape- and test-functions. While on the one hand this is not a real problem when the equation at hand is a simple one such as \eqref{Poisson}, on the other hand the choice of the right integration formula can turn to a sort of headache, in the case of more complex nonlinear PDEs.  \\        
\indent (vi) A priori the finite-element methodology studied in this article to solve boundary value problems posed in smooth curved domains is an advantageous alternative in many respects to more classical techniques such as parametric versions of the finite element method. This is because its most outstanding features are universality and simplicity, and eventually accuracy and CPU time too. We refer to the Appendix hereafter 
for some material related to the latter aspects. \\

To close this article, it is not superfluous to stress again that the technique advocated in this work to handle Dirichlet conditions prescribed on curvilinear boundaries has a wide scope of applicability. This is particularly true of three-dimensional boundary value problems. The author refers to \cite{CAMWA} for studies thereof applied to the equations of incompressible viscous flow, besides \cite{arXiv3D} and \cite{IMAJNA}. Applications to linear elasticity problems can be found in \cite{CFM2017} and \cite{Maugin}. The case of mixed finite element methods with normal-component degrees of freedom, such as Raviart-Thomas elements \cite{RaviartThomas}, will be addressed in a forthcoming paper. 

\section*{Acknowledgment}

The author gratefully acknowledges the financial support from CNPq, the National Research Council of Brazil, 
during the accomplishment of this work. 

\section*{APPENDIX - A comparison with the isoparametric technique}

Besides the advantages of the method studied in this work over the isoparametric technique (when it exists) 
pointed out in Section 5, the former appears to be more accurate than the latter. In support to this assertion,  we borrowed some numerical results from \cite{PAMM}, aimed at comparing our approach with the isoparametric version of the finite element method in terms of accuracy. \\
This comparison was carried out by solving with both methods for $k=2$ the convection-diffusion equation \eqref{Poisson} in the ellipse delimited by the curve of equation $(x/e)^2+y^2 = 1$, and for the exact solution $u=( e^2- e^2 x^2 - y^2)(e^2 - x^2 - e^2 y^2)$.  \\
We took $\nu=1$, ${\bf b} =(x,-y)$, $e=0.5$ and $g \equiv 0$, and owing to symmetry we considered only the quarter domain given by $x>0$ and $y>0$ by prescribing Neumann boundary conditions on $x=0$ and $y=0$.\\
In the table below the $L^2(\Omega_h)$-norms of both the gradient of the error function and 
of this function itself, and the pseudo-$L^{\infty}$-semi-norm $\parallel \cdot \parallel_{0,\infty,h}$ of the error are supplied. The latter is the maximum absolute value of the error at the mesh nodes. \\  
We computed with a family of quasi-uniform meshes defined by a single integer parameter $I$, containing $2I^2$ triangles, such that $h=1/I$. The notation $\tilde{u}_h$ is employed to represent the isoparametric solution. \\
As one can infer from these results the observed convergence rates are the same for both methods in the three 
senses under consideration, and furthermore conform to what can be expected from them. However the new method was more accurate than the isoparametric technique all the way. In particular, results for coarse meshes and approximations of the nodal values for all meshes are significantly better with the new method. \\
  
\begin{table*}[h]
{\small 
\centering
\begin{tabular}{ccccccc} &\\ [-.3cm]  
$h$ & $\longrightarrow$ & $1/4$ & $1/8$ & $1/16$ & $1/32$ & $1/64$ 
\tabularnewline & \\ [-.3cm] \hline &\\ [-.3cm]
$\parallel {\bf grad}(u-u_h) \parallel_{0,h}$ & $\longrightarrow$ & 0.539159 E-2 & 0.143611 E-2 & 0.367542 E-3 & 0.927845 E-4 & 0.233003 E-4   
\tabularnewline &\\ [-.3cm] \hline &\\ [-.3cm] 
$\parallel {\bf grad}(u-\tilde{u}_h) \parallel_{0,h}$ & $\longrightarrow$ &  0.793654 E-2 & 0.161494 E-2 & 0.392169 E-3 &  0.960439 E-4 & 0.237206 E-4 
\tabularnewline &\\ [-.3cm] \hline &\\ [-.3cm]
\tabularnewline &\\ [-.3cm] \hline &\\ [-.3cm]
$\parallel u-u_h \parallel_{0,h}$ & $\longrightarrow$ & 0.151255 E-3 & 0.184403 E-4 & 0.230467 E-5 & 0.289398 E-6 & 0.363189 E-7 
\tabularnewline &\\ [-.3cm] \hline &\\ [-.3cm] 
$\parallel u-\tilde{u}_h \parallel_{0,h}$ & $\longrightarrow$ & 0.173898 E-3 & 0.206144 E-4  & 0.247713 E-5 &  0.301543 E-6 & 0.371217 E-7
\tabularnewline &\\ [-.3cm] \hline &\\ [-.3cm]
\tabularnewline &\\ [-.3cm] \hline &\\ [-.3cm]
$\parallel u-u_h \parallel_{0,\infty,h}$ & $\longrightarrow$ & 0.397339 E-3 & 0.751885 E-4 & 0.110067 E-4 & 0.148037 E-5 & 0.195523 E-6  
\tabularnewline & \\ [-.3cm] \hline &\\ [-.3cm]
$\parallel u-\tilde{u}_h \parallel_{0,\infty,h}$ & $\longrightarrow$ & 0.753607 E-3 & 0.152953 E-3 & 0.241789 E-4 & 0.340427 E-5&0.452210 E-6
\tabularnewline &\\ [-.3cm] \hline &\\ [-.3cm]
\end{tabular} 
\caption{Errors with the new and the isoparametric approach for a problem in an ellipse and $k=2$.} 
}
\end{table*} 

Finally we report that both methods are roughly equivalent in terms of CPU time, as 
shown by a table supplied in \cite{arXiv2D}, for the same type of problem. 

\end{document}